\documentclass[11pt]{article}
\usepackage{latexsym}
\usepackage{bbm}
\usepackage{amscd}
\usepackage{amsmath}
\usepackage{amsfonts}
\usepackage{amssymb}
\usepackage{latexsym}
\usepackage[dvips]{graphicx}
\usepackage[colorlinks, citecolor=red]{hyperref}
\usepackage[marginal]{footmisc}
\usepackage[all]{xy}
\usepackage{verbatim}

\newtheorem{theorem}{\bf Theorem}[section]

\newtheorem{remark}[theorem]{\bf Remark}
\newenvironment{proof}{\noindent{}}{\quad \hfill$\Box$\vspace{2ex}}
\newenvironment{proof1}{\noindent{\em Proof:}}{\quad \hfill$\Box$\vspace{2ex}}

\def\no{\noindent}


\numberwithin{equation}{section}

\begin{document}

\thispagestyle{empty}
\begin{center}
{\LARGE \bf
Constantly curved holomorphic two-spheres
in the complex Grassmannian $G(2,6)$
with constant square norm of the second fundamental form
\no
}
\end{center}

\begin{center}
Jie Fei
\footnote{
J. Fei
\\School of Mathematics and Physics, Xi'an Jiaotong-Liverpool University, Suzhou 215123, P. R. China
\\E-mail: jie.fei@xjtlu.edu.cn
}
Ling He
\footnote{
L. He (Corresponding author)
\\Center for Applied Mathematics and KL-AAGDM, Tianjin University, Tianjin, 300072, P. R. China
\\E-mail: heling@tju.edu.cn
}
Jun Wang
\footnote{
J. Wang
\\School of Mathematics Sciences and Institute of Mathematics, Nanjing Normal University, Nanjing 210023, P. R. China
\\E-mail: wangjun706@njnu.edu.cn
}
\end{center}

\begin{center}
\parbox{12cm}
{\footnotesize{\bf ABSTRACT.}  We completely classify all noncongruent linearly full totally unramified constantly curved holomorphic two-spheres in $G(2,6)$ with constant square norm of the second fundamental form. They turn out to be homogeneous.
}
\end{center}

\no

\no

\section{Introduction}
\label{sec1}

Motivated by Calabi's rigidity principle about holomorphic isometric embedding from one complex manifold into the complex projective space in \cite{Calabi 1953}, one can study the holomorphic isometric embedding from one complex manifold into the complex Grassmannian. A holomorphic isometry starting from Riemann sphere
$S^2$ with metric of constant curvature is also called a constantly curved holomorphic two-sphere.
In differential geometry, the question of classification of noncongruent constantly curved holomorphic two-spheres in the complex Grassmannians is an important and difficult problem.
More interesting, it has close contact with the Grassmannian sigma models in theoretical physics (for example, see \cite{Delisle-Hussin-Zakrzewski 2013}).

We denote by $g_{FS}$ the Fubini-Study metric of constant holomorphic sectional curvature $4$.
From Calabi's rigidity principle (cf. \cite{Calabi 1953}; see also \cite{Rigoli 1985}, Sec.5 in \cite{Bolton-Jensen-Rigoli-Woodward 1988}), we know that if
$\varphi:S^2\rightarrow \left(\mathbb{C}P^{n},g_{FS}\right)$
is a linearly full (the image does not lie in some subspace $\mathbb{C}^{k}$ for $k<n+1$) constantly curved holomorphic immersion of degree $\mathrm{d}$, then $\mathrm{d}=n$ and $\varphi$ is the well-known Veronese embedding (up to $U(n+1)$) defined by
$$
V_0^{(n)}(z)=\left[1,\sqrt{{n\choose 1}}z,\cdots,\sqrt{{n\choose n}}z^n\right].
$$
It is interesting that $V_0^{(n)}(z)$ can generate the harmonic sequence
$V_0^{(n)},V_1^{(n)},\cdots,V_n^{(n)}$ in $\mathbb{C}P^n$,
where $V_i^{(n)}=[f_{i,0},\cdots,f_{i,p}\ ,\cdots,f_{i,n}]$ with $f_{i,p}$ being explicitly given by
$$f_{i,p}(z)=\frac{i!}{(1+z\bar{z})^i}\sqrt{\binom{n}{p}}z^{p-i}\sum_k(-1)^k \binom{p}{i-k}\binom{n-p}{k}(z\bar{z})^k.$$
Such a map $V_i^{(n)}:S^2\rightarrow \mathbb{C}P^n$ is a conformal minimal immersion with constant curvature and constant K\"ahler angle, which are given by
$$K_i^{(n)}=\frac{4}{n+2i(n-i)},\;\;\;\cos\alpha_i^{(n)}=\frac{n-2i}{n+2i(n-i)}.$$
This harmonic sequence is well known as  \textit{Veronese sequence} (cf. \cite{Bolton-Jensen-Rigoli-Woodward 1988}, Sec.5).

We denote by $ds^2_{G(m,m+n)}$ the standard K\"ahler metric of complex Grassmannian $G(m,m+n)~(m\leq n)$, which can also be induced from $g_{FS}$ by the Pl\"ucker embedding.
Chi-Zheng \cite{Chi-Zheng 1989} classified the constantly curved holomorphic two-spheres of degree $4$ in $\left(G(2,4),ds^2_{G(2,4)}\right)$
into two families (up to $U(4)$) by using the method of moving frames and Cartan's theory of higher order invariants (cf.\cite{Jensen 1977}).
This result means that the constantly curved holomorphic two-spheres in $\left(G(m,m+n),ds^2_{G(m,m+n)}\right)$ are more complicated.
Li-Yu \cite{Li-Yu 1999} (see \cite{Li-Jin 2008} for a detailed proof) showed that if $\varphi:S^2\rightarrow G(2,4)$
is a linearly full constantly curved holomorphic immersion of degree $\mathrm{d}$, then $\mathrm{d}=2,3,4$ and $\varphi$ is
explicitly characterized (up to $U(4)$).
For the case of $G(2,5)$, Jiao-Peng (\cite{Jiao-Peng 2004}, \cite{Jiao-Peng 2011})
proved that
if $\varphi:S^2\rightarrow G(2,5)$
is a non-singular constantly curved holomorphic immersion, then $\mathrm{d}=1,2,3,4,5$ and $\varphi$ is
explicitly characterized (up to $U(5)$) by using Pl\"ucker embedding.

Delisle-Hussin-Zakrzewski \cite{Delisle-Hussin-Zakrzewski 2013} recovered the known classification results in $G(2,4)$ and $G(2,5)$ mentioned above from the viewpoint of Grassmannian sigma models, and proposed two conjectures: Conjecture 1: \textit{If
$\varphi:S^2\rightarrow G(m,m+n)$
is a constantly curved holomorphic immersion of degree $\mathrm{d}$, then $\mathrm{d}\leq mn$.}
Conjecture 2: \textit{For $m,n$ fixed, the holomorphic immersion $\varphi$ can be constructed for $1\leq\mathrm{d}\leq mn$.}

Until now, there are few results about these two problems.
Under the assumption of homogeneity (the image is an orbit of an isometry subgroup of the target space), Peng-Xu \cite{Peng-Xu 2014} and Fei \cite{Fei 2019} independently used the representations of $SU(2)$ to give a complete classification of linearly full constantly curved holomorphic two-spheres of degree $\mathrm{d}$ in $G(2,n+2)$
and obtained that $\mathrm{d}$ takes $n$ or $2n$.
Delisle-Hussin-Zakrzewski (\cite{Delisle-Hussin-Zakrzewski 2013}, Proposition 1) declared that $\mathrm{d}\leq 6$ for $G(2,5)$. Recently the second named author gave some detailed discussions in \cite{He 2022}. Interestingly, Chi-Xie-Xu \cite{Chi-Xie-Xu 2022} constructed many non-homogeneous and singular (in Jiao-Peng's sense) constantly curved holomorphic two-spheres of degree $6$ in $G(2,5)$ and described the moduli space of such two-spheres.

S.S. Chern et al.(\cite{Chern 1968}, \cite{Chern-do Carmo-Kobayashi 1970}) showed the interest of the study of compact minimal submanifolds on the unit sphere $S^{n+p}$ with constant square norm of the second fundamental form, which led to the famous Chern conjecture for discreteness of such constant values.
N. Mok et al.(\cite{Mok 2009}, \cite{Mok-Ng 2009}) indicated that the second fundamental form can be used as a tool for studying the non-totally geodesic holomorphic isometric embeddings of the Poincar\'{e} disk into bounded symmetric domains.
Inspired by these works,
the second fundamental forms are expected to play an important role for studying constantly curved holomorphic two-spheres in $G(m,m+n)$.
Recently, the first and third named authors and Xu (\cite{WFX 2022},\cite{FW 2024})
applied the method of moving frames and harmonic sequences to completely classify
constantly curved holomorphic two-spheres in $G(2,N)$ and $G(3,N)$ with the square norm of the second fundamental form satisfying a certain pinching condition.

The current article is an attempt to study linearly full constantly curved holomorphic two-spheres of degree $\mathrm{d}$ in $G(2,n+2)$
with constant square norm of the second fundamental form.
For the reducible case, the first named author, Jiao and Xu in \cite{Fei-Jiao-Xu 2011} and the second named author in her Ph.D. Thesis \cite{He 2014}
showed that $\mathrm{d}=n$ with $\varphi=V_0^{(n+1)}\oplus V_1^{(n+1)}$ (up to $U(n+2)$) or $\mathrm{d}=2n$ with $\varphi=V_0^{(n)}\oplus v$ (up to $U(n+2)$), where $v$ is a non-zero constant vector.
Notice that there do not need any assumption with respect to the second fundamental form,
and these two holomorphic curves are both homogeneous.

For irreducible case, we have the Gauss equation (cf.\cite{Fei 2019},\cite{Fei-Xu 2017})
\begin{equation}\label{Gauss-eq0}
K=4-8|\det A_1|^2-\frac{S}{2},
\end{equation}
where $K={4}/{\mathrm{d}}$,
$S$ is the square norm of the second fundamental form,
and $|\det A_1|^2\phi^2\overline{\phi}^2$ is a global defined $(2,2)$-form (see section 2).
The key idea is to transform the condition of $S$ being constant into the one
that $|\det A_1|^2$ is a constant. Notice that $S$ involving the second derivative of the holomorphic mapping is hard to compute, while $|\det A_1|^2$ is only related to the first derivative of the holomorphic mapping. This makes it possible to calculate explicitly.
Thus by solving algebraic equations in the congruence class, we get the following main theorem.
\begin{theorem}\label{main-theorem1}
Let $\varphi:S^2\rightarrow G(2,n+2)$ be a linearly full irreducible constantly curved holomorphic two-sphere of degree $\mathrm{d}=4$ with constant square norm of the second fundamental form $S$.
Then $n=4$ and the corresponding holomorphic two-spheres are given as follows, up to $U(6)$,
\begin{equation}\label{eq0}
\varphi=\begin{bmatrix}
1&0&|t-2|z^3&\sqrt{t}z&\sqrt{\frac{3-t}{4-t}}|t-2|z^2
&\sqrt{\frac{t}{4-t}}z^2\\
0&1&\sqrt{3-t}z^2&0&\sqrt{4-t}z&0
\end{bmatrix}\subset G(2,6),
~0<t\leq 3
\end{equation}
with $S=t^2-4t+6$. Here when $t=3$, $\varphi=V_0^{(1)}\oplus V_0^{(3)}$ with $S=3$;
when $t=2$, $\varphi=V_0^{(2)}\oplus V_0^{(2)}$ with $S=2$.
\end{theorem}

Applying Theorem \ref{main-theorem1} and harmonic sequence, we obtain
the classification theorem of constantly curved holomorphic two-spheres in $G(2,6)$
with constant square norm of the second fundamental form, as follows.
\begin{theorem}\label{main-theorem2}
Let $\varphi:S^2\rightarrow G(2,6)$ be a linearly full constantly curved holomorphic two-sphere of degree $\mathrm{d}$ with constant square norm of the second fundamental form.
\\(1) If $\varphi$ is reducible, then $\mathrm{d}=4$ with $\varphi=V_0^{(5)}\oplus V_1^{(5)}$ (up to $U(6)$),
or $\mathrm{d}=8$ with $\varphi=V_0^{(4)}\oplus v$ (up to $U(6)$), where $v$ is a non-zero constant vector.
\\(2) If $\varphi$ is irreducible and totally unramified, then $\mathrm{d}=4$ and
$\varphi$ is unitary congruent to \eqref{eq0}.
\end{theorem}

\begin{remark}
It follows from Theorem \ref{main-theorem2}
that DHZ's conjecture 1 is true for $G(2,6)$
in the case of totally unramified and constant square norm of the second fundamental form.
\end{remark}

Recently, we \cite{Fei-He-Wang 2023} give a local rigidity characterization of all homogeneous holomorphic two-spheres in $G(2,N)$ in terms of a new global invariant $\kappa$ defined by the square norm of $(1,0)$ part of the second order covariant differential of the first $\partial$-transform for a holomorphic curve in $G(2,N)$.
Specifically, we showed that a linearly full irreducible constantly curved holomorphic two-sphere in $G(2,N)$ with constant square norm of the second fundamental form and $\kappa$ vanishing identically
is unitary congruent to
$$
V_0^{(n)}\oplus \left(\cos\theta V_1^{(n)}+\sin\theta V_0^{(n-2)}\right)
~\text{or}
~V_0^{(n_1)}\oplus V_0^{(n_2)},
$$
where $\theta\in(0,\frac{\pi}{2}]$,
$N=2n$ or $N=n_1+n_2+2$.
By virtue of the explicit expression given by \eqref{eq0} and straightforward calculations, we checked that the embedding \eqref{eq0} has vanishing $\kappa$ for each $t$, which concludes that the holomorphic curve given by \eqref{eq0} is unitary congruent to the homogeneous one
$$
V_0^{(3)}\oplus \left(\cos\theta V_1^{(3)}+\sin\theta V_0^{(1)}\right)
$$
and the correspondence between the two parameters is that $t=3\sin^2\theta$.
Therefore, the classification results of Theorem \ref{main-theorem2} tell us that
\begin{theorem}\label{main-theorem3}
A linearly full totally unramified constantly curved holomorphic two-sphere in $G(2,6)$
with constant square norm of the second fundamental form
must be homogeneous.
\end{theorem}

The requirement that a constantly curved holomorphic two-sphere is totally unramified is a somehow strong assumption.
Thus one is optimistic about generalizing Theorem \ref{main-theorem3} to any complex Grassmannian, which deserves further investigation.
In addition, whether is the assumption that the constantly curved holomorphic two-sphere is totally unramified in Theorem \ref{main-theorem2} and Theorem \ref{main-theorem3} necessary?
Since the condition of constant square norm of the second fundamental form leads to that
the constantly curved holomorphic two-sphere is unramified,
then can one construct an example of non-homogeneous constantly curved holomorphic two-sphere with constant square norm of the second fundamental form satisfying that the second element of the harmonic sequence has ramified points?

This paper is organized as follows.
In section 2, we give a characterization of constantly curved holomorphic two-sphere in $G(2,n+2)$
with constant square norm of the second fundamental form in terms of a system of algebraic equations. In section 3, we firstly solve the algebraic system to prove Theorem \ref{main-theorem1}.
Applying this result, together with the theory of harmonic sequence, we verify Theorem \ref{main-theorem2}.

{\bf{Acknowledgments}}~
This work is supported by National Key R\&D Program of China No. 2022YFA1006600 and NSF in China Nos. 12071352, 12071338, 11401481, 11301273, 11971237.
The first named author was also supported by the Research Enhancement Fund of Xi'an Jiaotong-Liverpool University (REF-18-01-03). The third named author was also supported by the NSF of the Jiangsu Higher Education Institutions of China (Grant No.~17KJA110002, No.~19KJA320001) and the Natural Science Foundation of Jiangsu Province (BK20181381).

\section{Constantly curved holomorphic two-spheres
in $G(2,n+2)$}

Let $\varphi:S^2\rightarrow G(2,n+2)$ be a linearly full constantly curved holomorphic two-sphere of degree $\mathrm{d}$.
Let $\varphi(0)=\mathcal{Z}_0\in G(2,n+2)$.
We choose a rectangular coordinate system in $\mathbb{C}^{n+2}$
such that the extended matrix $\mathcal{Z}_0=(I_2,\mathbf{0})$.
Then on the neighborhood of $\mathcal{Z}_0$ denoted by $\mathcal{V}_{\mathcal{Z}_0}$, we can write
$$
\varphi(z)=[I_2,F(z)],
$$
where $F(z)=\begin{pmatrix}
F_1(z)\\
F_2(z)
\end{pmatrix}$ is a $2\times n$ matrix-valued holomorphic function satisfying $F(0)=0$.

Let
$$Pl:\left(G(2,n+2),ds^2_{G(2,n+2)}\right)\rightarrow \left(\mathbb{C}P^N,g_{FS}\right)~(N=\frac{(n+2)(n+1)}{2}-1)
$$ be the standard Pl\"ucker embedding, which is a holomorphic isometric embedding (cf. \cite{Griffiths-Harris 1978}).
Then $Pl\circ \varphi:S^2\rightarrow \mathbb{C}P^N$
given by
$$
Pl\circ \varphi=[v_1\wedge v_2]=\begin{bmatrix}
1&F_2&-F_1&F_1\wedge F_2
\end{bmatrix}
$$
is a constantly curved holomorphic two-sphere, which may be not linearly full.
It follows from Calabi's rigidity theorem (cf. \cite{Calabi 1953}; see also \cite{Rigoli 1985}, Sec.5 in \cite{Bolton-Jensen-Rigoli-Woodward 1988}) that there exists a constant matrix $U\in U(N+1)$ such that
$Pl\circ \varphi=V_0^{(\mathrm{d})}\cdot U$,
where $V_0^{(\mathrm{d})}:S^2\rightarrow \mathbb{C}P^N$ is the holomorphic Veronese embedding given by
$$
V_0^{(\mathrm{d})}(z)=\begin{bmatrix}
1&\sqrt{\mathrm{d}\choose 1}z&\cdots&\sqrt{\mathrm{d}\choose k}z^k&\cdots &\sqrt{\mathrm{d}\choose \mathrm{d}}z^{\mathrm{d}} &0&\cdots &0
\end{bmatrix}.
$$
We immediately conclude that $\mathrm{d}\leq N$ and
\begin{equation}\label{cc}
1+|F_1|^2+|F_2|^2+|F_1\wedge F_2|^2=(1+z\bar{z})^{\mathrm{d}}.
\end{equation}
Moreover, we have the Gauss equation (cf. \cite{Fei 2019},\cite{Fei-Xu 2017})
\begin{equation}\label{Gauss-eq}
K=4-8|\det A_1|^2-\frac{S}{2},
\end{equation}
where $K={4}/{\mathrm{d}}$,
$S$ is the square norm of the second fundamental form,
and $|\det A_1|^2\phi^2\overline{\phi}^2$ is a global defined $(2,2)$-form,
the explicit explanation is as follows,
here $\phi$ is a local unitary coframe of $(1,0)$ type with respect to the induced metric $ds^2=\varphi^*(ds^2_{G(2,n+2)})$.

Let $\varphi=\begin{bmatrix}
I_2~F
\end{bmatrix},
~F(z)=\begin{bmatrix}
F_1(z)\\
F_2(z)
\end{bmatrix}$.
Set $v_1=\begin{bmatrix}
1&0&F_1
\end{bmatrix}$,
$v_2=\begin{bmatrix}
0&1&F_2
\end{bmatrix}$,
$e_1=\frac{v_1}{|v_1|}$,
$e_2=\frac{v_2-\frac{\left\langle v_2,v_1\right\rangle}{\left\langle v_1,v_1\right\rangle}v_1}
{\left|v_2-\frac{\left\langle v_2,v_1\right\rangle}{\left\langle v_1,v_1\right\rangle}v_1\right|}$.
Let $e_3,e_4$ be a unitary basis of the orthogonal projection plane $\varphi^{\perp}(\partial_z \varphi)$.
Then we have
\begin{eqnarray*}
&\partial_z e_1=a_{11}e_1+a_{12}e_2+a_{13}e_3+a_{14}e_4,\\
&\partial_z e_2=a_{21}e_1+a_{22}e_2+a_{23}e_3+a_{24}e_4.
\end{eqnarray*}
It follows that
$$
|\det A_1|^2\phi^2\overline{\phi}^2=|a_{13}a_{24}-a_{14}a_{23}|^2dz^2d\bar{z}^2.
$$
If $|\det A_1|^2\phi^2\overline{\phi}^2$ is identically equal to zero on $S^2$, then $\varphi$ is called \emph{reducible}.
If $|\det A_1|^2\phi^2\overline{\phi}^2$ is not identically equal to zero on $S^2$, then $\varphi$ is called \emph{irreducible}.
In the latter case, $|\det A_1|^2\phi^2\overline{\phi}^2$ has isolated zeros, which is called \emph{ramified points};
the order of zeros is called \emph{multiplicity of the corresponding ramified points}.
In particular, if $|\det A_1|^2\phi^2\overline{\phi}^2$ has no zeros, then $\varphi$ is called \emph{unramified}.
Similarly, we can give the definition of unramified minimal two-sphere.
If every element of the harmonic sequence generated by $\varphi$ is unramified, then $\varphi$ is called \emph{totally unramified} (see Definition 2.4 in \cite{He-Jiao-Zhou 2015}).

Now we give a characterization of constantly curved holomorphic two-spheres with constant square norm of the second fundamental form as follows:
\begin{theorem}\label{holos-thm2}
Let $\varphi:S^2\rightarrow G(2,n+2)$ be a linearly full constantly curved holomorphic two-sphere of degree $\mathrm{d}$ with constant square norm of the second fundamental form.
Locally, set
$$
\varphi=\begin{bmatrix}
I_2~F
\end{bmatrix},
~F(z)=\begin{bmatrix}
F_1(z)\\
F_2(z)
\end{bmatrix}
$$
where
$F(z)$ is a $2\times n$ matrix-valued holomorphic function satisfying $F(0)=0$.
Then
\begin{equation}\label{holos-thm2-eq1}
1+|F_1|^2+|F_2|^2+|F_1\wedge F_2|^2=(1+z\bar{z})^{\mathrm{d}},
\end{equation}
and
\begin{equation}\label{holos-thm2-eq2}
|\partial_zF_1\wedge\partial_zF_2|^2+
|\partial_zF_1\wedge\partial_zF_2\wedge F_1|^2+
|\partial_zF_1\wedge\partial_zF_2\wedge F_2|^2
+|\partial_zF_1\wedge\partial_zF_2\wedge F_1\wedge F_2|^2
=c(1+z\bar{z})^{2\mathrm{d}-4},
\end{equation}
where $c$ is a non-negative constant.
\end{theorem}
\begin{proof1}
If $S$ is constant, then $|\det A_1|^2$ is constant by \eqref{Gauss-eq}.
Here $|v_1\wedge v_2|^2=(1+z\bar{z})^{\mathrm{d}}$
and $e_1\wedge e_2=\frac{v_1\wedge v_2}{|v_1\wedge v_2|}$.
On the one hand,
\begin{eqnarray*}
\partial_z(e_1\wedge e_2)
&=&\partial_z e_1\wedge e_2+e_1\wedge\partial_z e_2\\
&=&(a_{11}+a_{22})e_1\wedge e_2-a_{13}e_2\wedge e_3-a_{14}e_2\wedge e_4
+a_{23}e_1\wedge e_3+a_{24}e_1\wedge e_4,
\end{eqnarray*}
which means
\begin{eqnarray*}
\partial_z(e_1\wedge e_2)\wedge \partial_z(e_1\wedge e_2)
=2(-a_{13}a_{24}+a_{14}a_{23})e_1\wedge e_2\wedge e_3\wedge e_4.
\end{eqnarray*}
On the other hand,
$$
\partial_z(e_1\wedge e_2)\wedge \partial_z(e_1\wedge e_2)
=\frac{-2v_1\wedge v_2\wedge\partial_zv_1\wedge \partial_zv_2}{|v_1\wedge v_2|^2}.
$$
Since the induced metric is $ds^2=\phi\overline{\phi}
=\frac{\mathrm{d}}{(1+z\bar{z})^2}dzd\bar{z}$,
then we have
$$
|\det A_1|^2=\frac{\left|v_1\wedge v_2\wedge\partial_zv_1\wedge \partial_zv_2\right|^2}
{{\mathrm{d}}^2(1+z\bar{z})^{2\mathrm{d}-4}}.
$$
Setting $|\det A_1|^2={c}/{{\mathrm{d}}^2}$, where $c$ is a non-negative constant,
we get
\begin{equation}\label{eq2+}
\left|v_1\wedge v_2\wedge\partial_zv_1\wedge \partial_zv_2\right|^2=c(1+z\bar{z})^{2\mathrm{d}-4}.
\end{equation}
A straightforward calculation verifies \eqref{holos-thm2-eq2}.
Furthermore we have
$$
S=8-\frac{16c}{{\mathrm{d}}^2}-\frac{8}{\mathrm{d}}.
$$
\end{proof1}

Let $F=\sum\limits_{\alpha=1}^NA_{\alpha}z^\alpha,
~A_\alpha=\begin{bmatrix}
a_1^{(\alpha)}
\\ a_2^{(\alpha)}
\end{bmatrix},
~W_\alpha=\begin{bmatrix}
a_1^{(\alpha)}
&a_2^{(\alpha)}
\end{bmatrix}$.
Then we have
$$
F_1=\sum\limits_{\alpha=1}^Na_1^{(\alpha)}z^\alpha,
~F_2=\sum\limits_{\alpha=1}^Na_2^{(\alpha)}z^\alpha.
$$
Set
$$
F_1\wedge F_2=\sum\limits_{j=2}^{2N}V_jz^j,
~\partial_zF_1\wedge\partial_zF_2=\sum\limits_{j=2}^{2N}R_jz^{j-2}.
$$
A straightforward calculation shows
$$
\partial_zF_1\wedge\partial_zF_2\wedge F_1
=\sum\limits_{j=2}^{2N}\sum\limits_{\alpha=1}^NR_j\wedge
a_1^{(\alpha)}z^{j+\alpha-2}
=\sum\limits_{p=4}^{3N}S_pz^{p-2},
$$
$$
\partial_zF_1\wedge\partial_zF_2\wedge F_2
=\sum\limits_{j=2}^{2N}\sum\limits_{\alpha=1}^NR_j\wedge
a_2^{(\alpha)}z^{j+\alpha-2}
=\sum\limits_{p=4}^{3N}T_pz^{p-2},
$$
$$
\partial_zF_1\wedge\partial_zF_2\wedge F_1\wedge F_2
=\sum\limits_{j,k=2}^{2N}R_j\wedge V_kz^{j+k-2}
=\sum\limits_{p=6}^{4N}X_pz^{p-2}.
$$
Set
$$
U=\begin{bmatrix}
  1 &\mathbf{0} &\mathbf{0} \\
  0 &W_1 & \mathbf{0}\\
  \vdots &\vdots &\vdots\\
  0 &W_{\mathrm{d}} &V_{\mathrm{d}}
\end{bmatrix},
~
Q=\begin{bmatrix}
  R_2 &\mathbf{0} &\mathbf{0} &\mathbf{0} \\
  R_3 &\mathbf{0} &\mathbf{0} &\mathbf{0}\\
  R_4 &S_4 &T_4 &\mathbf{0}\\
  R_5 &S_5 &T_5 &\mathbf{0}\\
  R_6 &S_6 &T_6 &X_6\\
  \vdots &\vdots &\vdots &\vdots\\
  R_{2\mathrm{d}-2} &S_{2\mathrm{d}-2} &T_{2\mathrm{d}-2} &X_{2\mathrm{d}-2}
\end{bmatrix},
$$
then we see that $U\in\mathbb{C}^{(\mathrm{d}+1)\times(N+1)}$ and
$Q\in\mathbb{C}^{(2\mathrm{d}-3)\times\left({n\choose 2}+2{n\choose 3}+{n\choose 4}\right)}$.

Applying Theorem \ref{holos-thm2}, we know that \eqref{holos-thm2-eq1} is equivalent to
\begin{equation}\label{holos-thm3-eq1}
UU^{*}=\Lambda_1,
\end{equation}
where $*$ denotes conjugate transpose and
$$
\Lambda_1=\text{diag}\left\{1,{\mathrm{d}\choose 1},\cdots,{\mathrm{d}\choose \mathrm{d}}\right\},
$$
and \eqref{holos-thm2-eq2} is equivalent to
\begin{equation}\label{holos-thm3-eq2}
QQ^{*}=\Lambda_2,
\end{equation}
where
$$
\Lambda_2=\text{diag}\left\{c,{2\mathrm{d}-4\choose 1}c,\cdots,{2\mathrm{d}-4\choose 2\mathrm{d}-4}c\right\}.
$$

In order to determine the linearly full constantly curved holomorphic two-spheres of degree $\mathrm{d}$ in $G(2,n+2)$ with constant square norm of the second fundamental form,
we need to solve \eqref{holos-thm3-eq1} and \eqref{holos-thm3-eq2}, modulo extrinsically the ambient unitary $U(n+2)$-congruence.
Considering the case of $\mathrm{d}=4$, we can get our main result.

\section{Proof of Theorems \ref{main-theorem1} and \ref{main-theorem2}}

\begin{proof}{\emph{Proof of Theorem \ref{main-theorem1}}}
Substituting $\mathrm{d}=4$ into \eqref{holos-thm3-eq1} and \eqref{holos-thm3-eq2} respectively yields
\begin{eqnarray}\label{holos-thm3-eq5-III}
U=\begin{bmatrix}
  1 &\mathbf{0} &\mathbf{0} \\
  0 &W_1 & \mathbf{0}\\
  0 &W_2 &V_2\\
  0 &W_3 &V_3\\
  0 &W_4 &V_4
\end{bmatrix},
~UU^*=\text{diag}\left\{1,4,6,4,1\right\}
\end{eqnarray}
and
\begin{eqnarray}\label{holos-thm3-eq6-III}
Q=\begin{bmatrix}
  R_2 &\mathbf{0} &\mathbf{0} &\mathbf{0} \\
  R_3 &\mathbf{0} &\mathbf{0} &\mathbf{0}\\
  R_4 &S_4 &T_4 &\mathbf{0}\\
  R_5 &S_5 &T_5 &\mathbf{0}\\
  R_6 &S_6 &T_6 &X_6\\
\end{bmatrix},
~QQ^*=\text{diag}\left\{c,4c,6c,4c,c\right\}.
\end{eqnarray}
In the following, we will discuss \eqref{holos-thm3-eq5-III} and \eqref{holos-thm3-eq6-III}
by case I: $W_4=0=W_3$,
case II: $W_4=0,W_3\neq 0$ and case III: $W_4\neq 0$ respectively.

\textbf{Case I:} $W_4=0=W_3$.
In this case, from \eqref{holos-thm3-eq5-III} and \eqref{holos-thm3-eq6-III}, we get $|V_3|^2=4$ and $|R_3|^2=4c$ respectively. Since $R_3=2V_3$, then we obtain $c=4$. Using \eqref{holos-thm3-eq5-III} and \eqref{holos-thm3-eq6-III} again,
we have $|W_1|^2=|R_2|^2=4$, which implies
\begin{equation}\label{I-eq1}
  |a_1^{(1)}|^2+|a_{2}^{(1)}|^2=4,
~|a_1^{(1)}\wedge a_2^{(1)}|^2=4.
\end{equation}
In the congruent class,
using the singular-value decomposition of complex matrix $A_1$, we can take
\begin{eqnarray*}
&a_1^{(1)}=\begin{pmatrix}
a_{11}^{(1)}&0&0&0&\cdots&0
\end{pmatrix},
a_2^{(1)}=\begin{pmatrix}
0&a_{22}^{(1)}&0&0&\cdots&0
\end{pmatrix},\\
&a_1^{(2)}=\begin{pmatrix}
a_{11}^{(2)}&a_{12}^{(2)}&a_{13}^{(2)}&a_{14}^{(2)}&\cdots&a_{1n}^{(2)}
\end{pmatrix},
a_2^{(2)}=\begin{pmatrix}
a_{21}^{(2)}&a_{22}^{(2)}&a_{23}^{(2)}&a_{24}^{(2)}
&\cdots&a_{2n}^{(2)}
\end{pmatrix},
\end{eqnarray*}
where $a_{11}^{(1)}\geq a_{22}^{(1)}>0$.
It follows from \eqref{I-eq1} that
\begin{equation}\label{I-eq2}
a_{11}^{(1)}=a_{22}^{(1)}=\sqrt{2}.
\end{equation}
From \eqref{I-eq2}, together with $\left\langle W_1,~W_2\right\rangle=0$ by \eqref{holos-thm3-eq5-III}, we obtain
\begin{equation}\label{I-eq3}
a_{11}^{(2)}+a_{22}^{(2)}=0.
\end{equation}
Since $V_3=a_1^{(1)}\wedge a_2^{(2)}+a_1^{(2)}\wedge a_2^{(1)}$, then $|V_3|^2=4$ gives us $|a_1^{(1)}\wedge a_2^{(2)}+a_1^{(2)}\wedge a_2^{(1)}|^2=4$. Combining this with \eqref{I-eq2} and \eqref{I-eq3}, we have
\begin{equation}\label{I-eq4}
|a_1^{(2)}|^2+|a_2^{(2)}|^2-2|a_{11}^{(2)}|^2
-|a_{21}^{(2)}|^2-|a_{12}^{(2)}|^2=2.
\end{equation}
It follows from $V_2=R_2$ that $|V_2|^2=4$.
Since $|W_2|^2+|V_2|^2=6$ by \eqref{holos-thm3-eq5-III},
then we get $|W_2|^2=2$,
i.e.,
\begin{equation}\label{I-eq5}
|a_{1}^{(2)}|^2+|a_{2}^{(2)}|^2=2.
\end{equation}
Substituting \eqref{I-eq5} into \eqref{I-eq4} yields
$$
2|a_{11}^{(2)}|^2
+|a_{21}^{(2)}|^2+|a_{12}^{(2)}|^2=0,
$$
which implies
\begin{equation}\label{I-eq6}
a_{11}^{(2)}=0=a_{21}^{(2)}=a_{12}^{(2)}
=a_{22}^{(2)}.
\end{equation}
From \eqref{holos-thm3-eq5-III}, we know $|V_4|^2=1$, i.e.,
$|a_1^{(2)}\wedge a_2^{(2)}|^2=1$. Combining this with \eqref{I-eq5} yields
\begin{equation}\label{I-eq7}
|a_{1}^{(2)}|^2=|a_{2}^{(2)}|^2=1,
~\left\langle a_{1}^{(2)},~a_{2}^{(2)}\right\rangle=0.
\end{equation}
Hence using \eqref{I-eq6} and \eqref{I-eq7}, we can take, in the congruent class,
\begin{eqnarray*}
&a_1^{(1)}=\begin{pmatrix}
\sqrt{2}&0&0&0&\cdots&0
\end{pmatrix},
a_2^{(1)}=\begin{pmatrix}
0&\sqrt{2}&0&0&\cdots&0
\end{pmatrix},\\
&a_1^{(2)}=\begin{pmatrix}
0&0&1&0&\cdots&0
\end{pmatrix},
a_2^{(2)}=\begin{pmatrix}
0&0&0&1
&\cdots&0
\end{pmatrix}.
\end{eqnarray*}
Then in this case, $\varphi$ is congruent to
$$
\varphi=\begin{bmatrix}
1&0&\sqrt{2}z&0&z^2&0\\
0&1&0&\sqrt{2}z&0&z^2
\end{bmatrix}
=V_0^{(2)}\oplus V_0^{(2)}\subset G(2,6).
$$

\textbf{Case II:} $W_4=0,W_3\neq 0$.
In this case, we have $V_6=a_1^{(3)}\wedge a_2^{(3)}=0$ and $V_5=a_1^{(2)}\wedge a_2^{(3)}+a_1^{(3)}\wedge a_2^{(2)}=0$. From \eqref{holos-thm3-eq5-III}, it follows that $\left\langle W_1,~W_3\right\rangle=0$.
Then in the congruent class,
we can take
\begin{eqnarray*}
&a_1^{(3)}=\begin{pmatrix}
a_{11}^{(3)}&0&0&0&\cdots&0
\end{pmatrix},
a_2^{(3)}=0,\\
&a_1^{(2)}=\begin{pmatrix}
a_{11}^{(2)}&a_{12}^{(2)}&a_{13}^{(2)}&a_{14}^{(2)}&\cdots&a_{1n}^{(2)}
\end{pmatrix},
a_2^{(2)}=\begin{pmatrix}
a_{21}^{(2)}&0&0&0&\cdots&0
\end{pmatrix},\\
&a_1^{(1)}=\begin{pmatrix}
0&a_{12}^{(1)}&a_{13}^{(1)}&a_{14}^{(1)}&\cdots&a_{1n}^{(1)}
\end{pmatrix},
a_2^{(1)}=\begin{pmatrix}
a_{21}^{(1)}&a_{22}^{(1)}&a_{23}^{(1)}&a_{24}^{(1)}
&\cdots&a_{2n}^{(1)}
\end{pmatrix},
\end{eqnarray*}
where $a_{11}^{(3)}>0$.
From \eqref{holos-thm3-eq5-III} and \eqref{holos-thm3-eq6-III},
we get $\left\langle V_2,~V_4\right\rangle=0$ and $\left\langle R_2,~R_4\right\rangle=0$ respectively. Since
$R_2=V_2=a_1^{(1)}\wedge a_2^{(1)}$,
$V_4=a_1^{(2)}\wedge a_2^{(2)}+a_1^{(3)}\wedge a_2^{(1)}$
and
$R_4=4a_1^{(2)}\wedge a_2^{(2)}+3a_1^{(3)}\wedge a_2^{(1)}$,
then we get
$$
\left\langle a_1^{(1)}\wedge a_2^{(1)}, a_1^{(3)}\wedge a_2^{(1)}\right\rangle=0,
$$
and
$$
\left\langle a_1^{(1)}\wedge a_2^{(1)}, a_1^{(2)}\wedge a_2^{(2)}\right\rangle=0,
$$
which implies by using $\left\langle a_1^{(1)},~a_1^{(3)}\right\rangle=0$
and $\left\langle a_1^{(1)},~a_2^{(2)}\right\rangle=0$ respectively,
\begin{equation}\label{II-eq1}
  \left\langle a_1^{(1)},~a_2^{(1)}\right\rangle a_{21}^{(1)}\overline{a_{11}^{(3)}}=0,
\end{equation}
and
\begin{equation}\label{II-eq2}
  \left\langle a_1^{(1)},~a_1^{(2)}\right\rangle a_{21}^{(1)}\overline{a_{21}^{(2)}}=0.
\end{equation}
In \eqref{holos-thm3-eq5-III}, we have $\left\langle W_1,~W_2\right\rangle=0$, i.e.,
\begin{equation}\label{II-eq3}
  \left\langle a_1^{(1)},~a_1^{(2)}\right\rangle+ a_{21}^{(1)}\overline{a_{21}^{(2)}}=0.
\end{equation}
It follows from \eqref{II-eq2} and \eqref{II-eq3} that
\begin{equation}\label{II-eq4}
  \left\langle a_1^{(1)},~a_1^{(2)}\right\rangle=0
=a_{21}^{(1)}\overline{a_{21}^{(2)}}.
\end{equation}

\textbf{Claim:} $a_{21}^{(1)}=0$.

Otherwise if $a_{21}^{(1)}\neq0$, then $\left\langle a_1^{(1)},~a_2^{(1)}\right\rangle=0$ by \eqref{II-eq1} and
$a_{21}^{(2)}=0$ by \eqref{II-eq4}. The latter tells us $a_2^{(2)}=0$.
In \eqref{holos-thm3-eq5-III}, we have $|V_4|^2=1$.
Since $V_4=a_1^{(3)}\wedge a_2^{(1)}$, then
\begin{equation}\label{II-eq5}
  |a_1^{(3)}\wedge a_2^{(1)}|^2=1.
\end{equation}
In \eqref{holos-thm3-eq6-III}, we have $|R_5|^2+|S_5|^2+|T_5|^2=4c$.
Since  $R_5=T_5=0,~S_5=-2a_1^{(1)}\wedge a_2^{(1)}\wedge a_1^{(3)}$, then
\begin{equation}\label{II-eq6}
  |a_1^{(1)}|^2|a_1^{(3)}\wedge a_2^{(1)}|^2=c.
\end{equation}
It follows from \eqref{II-eq5} and \eqref{II-eq6} that
\begin{equation}\label{II-eq7}
  |a_1^{(1)}|^2=c.
\end{equation}
In \eqref{holos-thm3-eq5-III} and \eqref{holos-thm3-eq6-III}, we have
$|W_1|^2=4$ and $|R_2|^2=c$ respectively, i.e.,
\begin{equation}\label{II-eq8}
  |a_1^{(1)}|^2+|a_2^{(1)}|^2=4,
  ~|a_1^{(1)}|^2|a_2^{(1)}|^2=c.
\end{equation}
Combining \eqref{II-eq7} and \eqref{II-eq8}, we obtain
$|a_2^{(1)}|^2=1$ and $|a_1^{(1)}|^2=c=3$.
From \eqref{holos-thm3-eq5-III} and \eqref{holos-thm3-eq6-III}, we obtain
$|W_3|^2+|V_3|^2=4$ and $|R_3|^2=4c=12$ respectively.
Using these, together with $R_3=2V_3$, we get $|W_3|^2=1$,
i.e., $|a_1^{(3)}|^2=1$. Substituting this and $|a_2^{(1)}|^2=1$ into \eqref{II-eq5} yields
$\left\langle a_1^{(3)},~a_2^{(1)}\right\rangle=0$, which shows
$a_{11}^{(3)}\overline{a_{21}^{(1)}}=0$. It follows that $a_{21}^{(1)}=0$, which contradicts our assumption.
Hence we verify $a_{21}^{(1)}=0$.

Now we know $\left\langle a_1^{(3)},~a_2^{(1)}\right\rangle=\left\langle a_2^{(2)},~a_2^{(1)}\right\rangle=0$. The latter gives us $\left\langle a_1^{(2)},~a_1^{(1)}\right\rangle=0$ by $\left\langle W_2,~W_1\right\rangle=0$
in \eqref{holos-thm3-eq5-III}.

From \eqref{holos-thm3-eq5-III} and \eqref{holos-thm3-eq6-III},
we get $\left\langle V_3,~V_4\right\rangle=0$ and $\left\langle R_3,~R_4\right\rangle=0$ respectively. Since
$R_3=2V_3=2a_1^{(1)}\wedge a_2^{(2)}+2a_1^{(2)}\wedge a_2^{(1)}$,
$V_4=a_1^{(2)}\wedge a_2^{(2)}+a_1^{(3)}\wedge a_2^{(1)}$
and
$R_4=4a_1^{(2)}\wedge a_2^{(2)}+3a_1^{(3)}\wedge a_2^{(1)}$,
then we get
$$
\left\langle a_1^{(1)}\wedge a_2^{(2)}+a_1^{(2)}\wedge a_2^{(1)}, a_1^{(2)}\wedge a_2^{(2)}\right\rangle=0,
$$
and
$$
\left\langle a_1^{(1)}\wedge a_2^{(2)}+a_1^{(2)}\wedge a_2^{(1)}, a_1^{(3)}\wedge a_2^{(1)}\right\rangle=0,
$$
which implies, respectively,
\begin{equation}\label{II-eq1+}
  \left\langle a_2^{(1)},~a_1^{(2)}\right\rangle a_{11}^{(2)}\overline{a_{21}^{(2)}}=0,
\end{equation}
and
\begin{equation}\label{II-eq2+}
  |a_2^{(1)}|^2a_{11}^{(2)}\overline{a_{11}^{(3)}}-\left\langle a_1^{(1)},~a_2^{(1)}\right\rangle a_{21}^{(2)}\overline{a_{11}^{(3)}}=0.
\end{equation}
In the following, we will discuss this case from two subcases of $a_{21}^{(2)}=0$ and $a_{21}^{(2)}\neq0$ respectively.

\textbf{Subcase II1:} If $a_{21}^{(2)}=0$, i.e., $a_{2}^{(2)}=0$, then
since $a_2^{(1)}$ is non-zero due to that $\varphi$ is unramified, we obtain $a_{11}^{(2)}=0$ by using \eqref{II-eq2+}.
It follows from the proof of the above {Claim}, we know that
\eqref{II-eq5}, \eqref{II-eq6} and \eqref{II-eq8}
become, respectively,
\begin{equation}\label{II-eq5+}
  |a_1^{(3)}|^2|a_2^{(1)}|^2=1.
\end{equation}
\begin{equation}\label{II-eq6+}
  |a_1^{(1)}\wedge a_2^{(1)}|^2|a_1^{(3)}|^2=c.
\end{equation}
\begin{equation}\label{II-eq8+}
  |a_1^{(1)}|^2+|a_2^{(1)}|^2=4,
  ~|a_1^{(1)}\wedge a_2^{(1)}|^2=c.
\end{equation}
In light of \eqref{II-eq5+}-\eqref{II-eq8+}, we assert immediately that
$|a_{11}^{(3)}|^2=|a_2^{(1)}|^2=1$ and $|a_1^{(1)}|^2=3$.
On the one hand, since $|W_3|^2+|V_3|^2=4$, then we get by $|W_3|^2=|a_{11}^{(3)}|^2=1$ that
$|V_3|^2=3$, which implies $|R_3|^2=4|V_3|^2=12$. In terms of this, together with $|R_3|^2=4c$, we get $c=3$.
Substituting it into \eqref{II-eq8+} yields $|a_1^{(1)}\wedge a_2^{(1)}|^2=3$.
Combining this with $|a_1^{(1)}|^2=3$ and $|a_2^{(1)}|^2=1$, we obtain
$\left\langle a_1^{(1)},~a_2^{(1)}\right\rangle=0$.

On the other hand, since $|W_2|^2+|V_2|^2=6$ and $|V_2|^2=3$, then we have $|W_2|^2=3$, i.e., $|a_{1}^{(2)}|^2=3$. In light of this, together with
$|V_3|^2=|a_1^{(2)}\wedge a_2^{(1)}|^2=3$ and $|a_2^{(1)}|^2=1$, we obtain that $\left\langle a_1^{(2)},~a_2^{(1)}\right\rangle=0$.
Hence in the congruent class,
we can take
\begin{eqnarray*}
&a_1^{(3)}=\begin{pmatrix}
1&0&0&0&\cdots&0
\end{pmatrix},
a_2^{(3)}=0,\\
&a_1^{(2)}=\begin{pmatrix}
0&\sqrt{3}&0&0&\cdots&0
\end{pmatrix},
a_2^{(2)}=0,\\
&a_1^{(1)}=\begin{pmatrix}
0&0&\sqrt{3}&0&\cdots&0
\end{pmatrix},
a_2^{(1)}=\begin{pmatrix}
0&0&0&1
&\cdots&0
\end{pmatrix}.
\end{eqnarray*}
Then in this subcase, $\varphi$ is congruent to
$$
\varphi=\begin{bmatrix}
1&0&z^3&\sqrt{3}z^2&\sqrt{3}z&0\\
0&1&0&0&0&z
\end{bmatrix}
=V_0^{(1)}\oplus V_0^{(3)}\subset G(2,6).
$$

\textbf{Subcase II2:} If $a_{21}^{(2)}\neq 0$, then we first claim $a_{11}^{(2)}=0$. Otherwise if $a_{11}^{(2)}\neq0$, then \eqref{II-eq1+} leads to $\left\langle a_2^{(1)}, a_1^{(2)}\right\rangle=0$. Combining this with $\left\langle a_1^{(1)}, a_1^{(2)}\right\rangle=\left\langle a_1^{(1)}, a_2^{(2)}\right\rangle=\left\langle a_2^{(1)}, a_2^{(2)}\right\rangle=0$, we
find $\left\langle V_2, V_3\right\rangle=0$, which implies $\left\langle W_2, W_3\right\rangle=0$ by applying $\left\langle W_2, W_3\right\rangle+\left\langle V_2, V_3\right\rangle=0$ in \eqref{holos-thm3-eq5-III}.
It follows that $a_{11}^{(2)}=0$, which is a contradiction. Thus we verify $a_{11}^{(2)}=0$.

Now \eqref{II-eq2+} yields $\left\langle a_1^{(1)}, a_2^{(1)}\right\rangle=0$.
In \eqref{holos-thm3-eq5-III}, we have $|W_1|^2=4$ and $|W_3|^2+|V_3|^2=4$, i.e.,
\begin{equation}\label{II2-eq1}
  |a_1^{(1)}|^2+|a_{2}^{(1)}|^2=4
\end{equation}
and
\begin{equation}\label{II2-eq2}
  |a_{11}^{(3)}|^2+|a_{1}^{(1)}|^2
|a_{21}^{(2)}|^2
+|a_{1}^{(2)}\wedge a_{2}^{(1)}|^2=4.
\end{equation}
In \eqref{holos-thm3-eq6-III}, we know $|R_2|^2=c$ and $|R_6|^2+|S_6|^2+|T_6|^2+|X_6|^2=c$. Since  $R_6=T_6=0,~S_6=-a_1^{(2)}\wedge a_2^{(1)}\wedge a_1^{(3)}$
and
$X_6=a_1^{(1)}\wedge a_2^{(1)}\wedge a_1^{(2)}\wedge a_2^{(2)}$,
then
\begin{equation}\label{II2-eq3}
  |a_{1}^{(1)}|^2|a_{2}^{(1)}|^2=c
\end{equation}
and
\begin{equation}\label{II2-eq4}
  |a_{1}^{(2)}\wedge a_{2}^{(1)}|^2\left(|a_{11}^{(3)}|^2
+|a_{1}^{(1)}|^2|a_{21}^{(2)}|^2\right)=c.
\end{equation}
Observing \eqref{II2-eq1}-\eqref{II2-eq4}, we find that $|a_{1}^{(1)}|^2,|a_{2}^{(1)}|^2$ and $|a_{11}^{(3)}|^2+|a_{1}^{(1)}|^2|a_{21}^{(2)}|^2,
|a_{2}^{(1)}\wedge a_{1}^{(2)}|^2$ are respectively the roots of equation $x^2-4x+c=0$. Then the following two cases will happen:
(i) $|a_{1}^{(1)}|^2
=|a_{11}^{(3)}|^2+|a_{1}^{(1)}|^2|a_{21}^{(2)}|^2,
|a_{2}^{(1)}|^2=|a_{2}^{(1)}\wedge a_{1}^{(2)}|^2,$
or (ii) $|a_{1}^{(1)}|^2=|a_{2}^{(1)}\wedge a_{1}^{(2)}|^2,
|a_{2}^{(1)}|^2=|a_{11}^{(3)}|^2+|a_{1}^{(1)}|^2|a_{21}^{(2)}|^2$.

Set $|a_{1}^{(1)}|^2=t (0<t<4)$, then $|a_{2}^{(1)}|^2=4-t=\frac{c}{t}$, i.e., $c=4t-t^2$.
It follows from $|W_3|^2+|V_3|^2=4$ and $|R_3|^2=4|V_3|^2=4c$ that $|W_3|^2=4-c$, i.e., $|a_{11}^{(3)}|^2=4-c=(t-2)^2$.
If (i) happens, then $|a_{21}^{(2)}|^2=1-\frac{4-c}{t}$.
Using $|R_5|^2+|S_5|^2+|T_5|^2=4c$ in \eqref{holos-thm3-eq6-III}
and $R_5=0,~S_5=-2a_1^{(1)}\wedge a_2^{(1)}\wedge a_1^{(3)}
+2a_1^{(1)}\wedge a_1^{(2)}\wedge a_2^{(2)},
T_5=2a_2^{(1)}\wedge a_1^{(2)}\wedge a_2^{(2)}$,
we get $|a_{1}^{(1)}|^2|V_4|^2+|a_{2}^{(1)}\wedge a_{1}^{(2)}|^2|a_{21}^{(2)}|^2=c$.
Combining this with $|V_4|^2=1$ by
\eqref{holos-thm3-eq5-III},
we obtain $|a_{1}^{(1)}|^2+|a_{2}^{(1)}|^2|a_{21}^{(2)}|^2=c$.
Then
$$
t+(4-t)\left(1-\frac{4-c}{t}\right)=c,
$$
which implies $t=2$. It follows that $a_{11}^{(3)}=0$, which contradicts $W_3\neq 0$. So this case of (i) will not happen.

If (ii) happens, then we get
\begin{equation}\label{II2-eq5}
|a_{21}^{(2)}|^2=\frac{4-t}{t}-\frac{4-c}{t}=\frac{c}{t}-1=3-t,
0<t<3,t\neq 2.
\end{equation}
Using $|W_2|^2+|V_2|^2=6$ in \eqref{holos-thm3-eq5-III},
we have
$$
|a_{1}^{(2)}|^2+|a_{21}^{(2)}|^2
+|a_{1}^{(1)}|^2|a_{2}^{(1)}|^2=6,
$$
which gives us by \eqref{II2-eq3} and \eqref{II2-eq5}
\begin{equation}\label{II2-eq6}
|a_{1}^{(2)}|^2=6-c-(3-t)=t^2-3t+3.
\end{equation}
Since $\left\langle a_1^{(1)},~a_1^{(2)}\right\rangle=0=\left\langle a_1^{(1)},~a_2^{(1)}\right\rangle$,
then
in the congruent class,
we can take
\begin{eqnarray*}
&a_1^{(3)}=\begin{pmatrix}
|t-2|&0&0&0&0&\cdots&0
\end{pmatrix},
a_2^{(3)}=0,\\
&a_1^{(2)}=\begin{pmatrix}
0&0&a_{13}^{(2)}&a_{14}^{(2)}&0&\cdots&0
\end{pmatrix},
a_2^{(2)}=\begin{pmatrix}
\sqrt{3-t}&0&0&0&0&\cdots&0
\end{pmatrix},\\
&a_1^{(1)}=\begin{pmatrix}
0&\sqrt{t}&0&0&0&\cdots&0
\end{pmatrix},
a_2^{(1)}=\begin{pmatrix}
0&0&\sqrt{4-t}&0&0
&\cdots&0
\end{pmatrix},
\end{eqnarray*}
where $a_{13}^{(2)}\in\mathbb{C},a_{14}^{(2)}>0$.
It follows from $|a_{1}^{(1)}|^2=|a_{2}^{(1)}\wedge a_{1}^{(2)}|^2$ that
$t=|a_{14}^{(2)}|^2|a_{23}^{(1)}|^2$, which implies
\begin{equation}\label{II2-eq7}
a_{14}^{(2)}=\sqrt{\frac{t}{4-t}}.
\end{equation}
Then \eqref{II2-eq6} and \eqref{II2-eq7} tell us
\begin{equation}\label{II2-eq8}
|a_{13}^{(2)}|^2=|a_{1}^{(2)}|^2-|a_{14}^{(2)}|^2
=t^2-3t+3-\frac{t}{4-t}=\frac{(3-t)(2-t)^2}{4-t}.
\end{equation}

On the other hand, substituting $V_4=a_1^{(2)}\wedge a_2^{(2)}+a_1^{(3)}\wedge a_2^{(1)}$ into $|V_4|^2=1$ leads to
$$
|a_{1}^{(2)}|^2|a_{21}^{(2)}|^2
+|a_{11}^{(3)}|^2|a_{2}^{(1)}|^2-a_{13}^{(2)}\overline{a_{23}^{(1)}}
\overline{a_{11}^{(3)}}
a_{21}^{(2)}-\overline{a_{13}^{(2)}}a_{23}^{(1)}
a_{11}^{(3)}\overline{a_{21}^{(2)}}=1,
$$
which yields by a straightforward calculation
\begin{equation}\label{II2-eq9}
a_{13}^{(2)}+\overline{a_{13}^{(2)}}=2\sqrt{\frac{3-t}{4-t}}|t-2|.
\end{equation}
Comparing \eqref{II2-eq8} and \eqref{II2-eq9}, we find
$$
\text{Re} ~a_{13}^{(2)}=\sqrt{\frac{3-t}{4-t}}|t-2|=|a_{13}^{(2)}|,
$$
which implies $a_{13}^{(2)}\in\mathbb{R}$ and
$$
a_{13}^{(2)}=\sqrt{\frac{3-t}{4-t}}|t-2|.
$$
Thus in this subcase, $\varphi$ is congruent to
$$
\varphi=\begin{bmatrix}
1&0&|t-2|z^3&\sqrt{t}z&\sqrt{\frac{3-t}{4-t}}|t-2|z^2
&\sqrt{\frac{t}{4-t}}z^2\\
0&1&\sqrt{3-t}z^2&0&\sqrt{4-t}z&0
\end{bmatrix}\subset G(2,6),
~0<t<3, t\neq 2.
$$

\textbf{Case III:} $W_4\neq 0$.
We will prove that there does not exist this case by contradiction.
Suppose that this case happens, then
we have $V_8=a_1^{(4)}\wedge a_2^{(4)}=0$ and $V_7=a_1^{(3)}\wedge a_2^{(4)}+a_1^{(4)}\wedge a_2^{(3)}=0$. From \eqref{holos-thm3-eq5-III}, it follows that $\left\langle W_1,~W_4\right\rangle=0$.
Then in the congruent class,
we can take
\begin{eqnarray*}
&a_1^{(4)}=\begin{pmatrix}
a_{11}^{(4)}&0&0&0&\cdots&0
\end{pmatrix},
a_2^{(4)}=0,\\
&a_1^{(3)}=\begin{pmatrix}
a_{11}^{(3)}&a_{12}^{(3)}&a_{13}^{(3)}&a_{14}^{(3)}
&\cdots&a_{1n}^{(3)}
\end{pmatrix},
a_2^{(3)}=\begin{pmatrix}
a_{21}^{(3)}&0&0&0&\cdots&0
\end{pmatrix},\\
&a_1^{(2)}=\begin{pmatrix}
a_{11}^{(2)}&a_{12}^{(2)}&a_{13}^{(2)}&a_{14}^{(2)}
&\cdots&a_{1n}^{(2)}
\end{pmatrix},
a_2^{(2)}=\begin{pmatrix}
a_{21}^{(2)}&a_{22}^{(2)}&a_{23}^{(2)}&a_{24}^{(2)}&
\cdots&a_{2n}^{(2)}
\end{pmatrix},\\
&a_1^{(1)}=\begin{pmatrix}
0&a_{12}^{(1)}&a_{13}^{(1)}&a_{14}^{(1)}&\cdots&a_{1n}^{(1)}
\end{pmatrix},
a_2^{(1)}=\begin{pmatrix}
a_{21}^{(1)}&a_{22}^{(1)}&a_{23}^{(1)}&a_{24}^{(1)}
&\cdots&a_{2n}^{(1)}
\end{pmatrix},
\end{eqnarray*}
where $a_{11}^{(4)}>0$.
Since $V_5=a_1^{(2)}\wedge a_2^{(3)}
+a_1^{(3)}\wedge a_2^{(2)}+a_1^{(4)}\wedge a_2^{(1)}=0$,
then we get $R_5=6a_1^{(2)}\wedge a_2^{(3)}
+6a_1^{(3)}\wedge a_2^{(2)}+4a_1^{(4)}\wedge a_2^{(1)}
=-2a_1^{(4)}\wedge a_2^{(1)}$.
Substituting this and $R_2=a_1^{(1)}\wedge a_2^{(1)}$ into
$\left\langle R_2,~R_5\right\rangle=0$ in \eqref{holos-thm3-eq6-III}, we obtain
$$
\left\langle a_1^{(1)}\wedge a_2^{(1)}, a_2^{(1)}\wedge a_1^{(4)}\right\rangle=0,
$$
which implies by using $\left\langle a_1^{(1)}, a_1^{(4)}\right\rangle=0$
that
\begin{equation}\label{III-eq1}
 \left\langle a_1^{(1)},~a_2^{(1)}\right\rangle a_{21}^{(1)}\overline{a_{11}^{(4)}}=0.
\end{equation}
Similarly from $V_6=a_1^{(3)}\wedge a_2^{(3)}
+a_1^{(4)}\wedge a_2^{(2)}=0$,
we know $R_6=9a_1^{(3)}\wedge a_2^{(3)}
+8a_1^{(4)}\wedge a_2^{(2)}
=a_1^{(3)}\wedge a_2^{(3)}$.
Substituting this and $R_2=a_1^{(1)}\wedge a_2^{(1)}$ into
$\left\langle R_2,~R_6\right\rangle=0$ in \eqref{holos-thm3-eq6-III}, we get
$$
\left\langle a_1^{(1)}\wedge a_2^{(1)}, a_1^{(3)}\wedge a_2^{(3)}\right\rangle=0,
$$
which implies by $\left\langle a_1^{(1)}, a_2^{(3)}\right\rangle=0$
that
\begin{equation}\label{III-eq2}
  \left\langle a_1^{(1)},~a_1^{(3)}\right\rangle a_{21}^{(1)}\overline{a_{21}^{(3)}}=0.
\end{equation}
In \eqref{holos-thm3-eq6-III}, we have $\left\langle W_1,~W_3\right\rangle=0$, i.e.,
\begin{equation}\label{III-eq3}
  \left\langle a_1^{(1)},~a_1^{(3)}\right\rangle+ a_{21}^{(1)}\overline{a_{21}^{(3)}}=0.
\end{equation}
It follows from \eqref{III-eq2} and \eqref{III-eq3} that
\begin{equation}\label{III-eq4}
\left\langle a_1^{(1)},~a_1^{(3)}\right\rangle=0
=a_{21}^{(1)}\overline{a_{21}^{(3)}}.
\end{equation}

\textbf{Claim 1:} $a_{21}^{(1)}=0$.

Otherwise if $a_{21}^{(1)}\neq0$, then $\left\langle a_1^{(1)},~a_2^{(1)}\right\rangle=0$ by \eqref{III-eq1}
and $a_{21}^{(3)}=0$ by \eqref{III-eq4}. The latter tells us $a_2^{(3)}=0$.
Then by using $V_6=a_1^{(4)}\wedge a_2^{(2)}=0$ we have $a_2^{(2)}=a_{21}^{(2)}/a_{11}^{(4)}\cdot a_1^{(4)}$.
This gives us $\left\langle a_1^{(1)},~a_2^{(2)}\right\rangle=0$.
At this time, since $R_2=a_1^{(1)}\wedge a_2^{(1)}$ and $R_4=4a_1^{(2)}\wedge a_2^{(2)}+3a_1^{(3)}\wedge a_2^{(1)}$, then $\left\langle R_2,~R_4\right\rangle=0$ in \eqref{holos-thm3-eq6-III} gives
\begin{equation}\label{III-eq5}
4\left\langle a_1^{(1)}\wedge a_2^{(1)}, a_1^{(2)}\wedge a_2^{(2)}\right\rangle
+3\left\langle a_1^{(1)}\wedge a_2^{(1)}, a_1^{(3)}\wedge a_2^{(1)}\right\rangle=0.
\end{equation}
It follows from $\left\langle a_1^{(1)},~a_1^{(3)}\right\rangle=\left\langle a_1^{(1)},~a_2^{(1)}\right\rangle=0$ that
\begin{equation}\label{III-eq6}
\left\langle a_1^{(1)}\wedge a_2^{(1)}, a_1^{(3)}\wedge a_2^{(1)}\right\rangle=0.
\end{equation}
Substituting \eqref{III-eq6} into \eqref{III-eq5} yields
\begin{equation*}
\left\langle a_1^{(1)}\wedge a_2^{(1)}, a_1^{(2)}\wedge a_2^{(2)}\right\rangle
=0,
\end{equation*}
which implies by $\left\langle a_1^{(1)},~a_2^{(2)}\right\rangle=0$, that
\begin{equation}\label{III-eq7}
\left\langle a_1^{(1)},~a_1^{(2)}\right\rangle a_{21}^{(1)}\overline{a_{21}^{(2)}}=0.
\end{equation}
On the other hand, in \eqref{holos-thm3-eq5-III}, we have $\left\langle W_1,~W_2\right\rangle=0$, i.e.,
\begin{equation}\label{III-eq8}
\left\langle a_1^{(1)},~a_1^{(2)}\right\rangle+ a_{21}^{(1)}\overline{a_{21}^{(2)}}=0.
\end{equation}
It follows from \eqref{III-eq7} and \eqref{III-eq8} that
\begin{equation}\label{III-eq9}
\left\langle a_1^{(1)},~a_1^{(2)}\right\rangle=0
=a_{21}^{(1)}\overline{a_{21}^{(2)}},
\end{equation}
which means $a_{21}^{(2)}=0$ by $a_{21}^{(1)}\neq 0$, i.e., $a_2^{(2)}=0$.
Now by using $V_5=a_1^{(4)}\wedge a_2^{(1)}=0$, we have $a_2^{(1)}=a_{21}^{(1)}/a_{11}^{(4)}\cdot a_1^{(4)}$.
And $V_2,V_3,V_4$ reduce to, respectively,
$$
V_2=a_1^{(1)}\wedge a_2^{(1)},
~V_3=a_1^{(2)}\wedge a_2^{(1)},
~V_4=a_1^{(3)}\wedge a_2^{(1)}.
$$
Observing that $\left\langle V_2,~V_3\right\rangle=\left\langle V_2,~V_4\right\rangle=0$. Substituting this into
$\left\langle W_2,~W_3\right\rangle+\left\langle V_2,~V_3\right\rangle=0$
and $\left\langle W_2,~W_4\right\rangle+\left\langle V_2,~V_4\right\rangle=0$
in \eqref{holos-thm3-eq5-III} respectively yields
$
\left\langle W_2,~W_3\right\rangle
=\left\langle W_2,~W_4\right\rangle
=0,
$
i.e.,
$
a_{11}^{(2)}=0,
~\left\langle a_1^{(2)},~a_1^{(3)}\right\rangle
=0.
$
This just gives us $\left\langle V_3,~V_4\right\rangle=0$.
Substituting it into
$\left\langle W_3,~W_4\right\rangle+\left\langle V_3,~V_4\right\rangle=0$
in \eqref{holos-thm3-eq5-III} yields
$\left\langle W_3,~W_4\right\rangle
=0,$ i.e., $a_{11}^{(3)}=0.$

Since $R_5=T_5=0$ and $S_5=-2a_1^{(1)}\wedge a_2^{(1)}\wedge a_1^{(3)}$,
then $|R_5|^2+|S_5|^2+|T_5|^2=4c$ in \eqref{holos-thm3-eq6-III}
leads to $|R_2|^2|a_1^{(3)}|^2=c$, which implies $|a_1^{(3)}|^2=1$
by $|R_2|^2=c$ in \eqref{holos-thm3-eq6-III}.
This tells us $|W_3|^2=1$. Substituting it into $|W_3|^2+|V_3|^2=4$ in \eqref{holos-thm3-eq5-III} yields $|V_3|^2=3$.
On the other hand, it follows from $|R_3|^2=4|V_3|^2=4c$ in \eqref{holos-thm3-eq6-III} that $|V_3|^2=c$.
Hence we get $c=3$, which implies $|V_2|^2=3$. Substituting it into $|W_2|^2+|V_2|^2=6$ in \eqref{holos-thm3-eq5-III} yields $|W_2|^2=3$,
i.e., $|a_1^{(2)}|^2=3$. Using this and $|V_3|^2=|a_1^{(2)}|^2|a_2^{(1)}|^2=3$ leads to $|a_2^{(1)}|^2=1$. In terms of this, together with $|a_1^{(3)}|^2=1$, we obtain $|V_4|^2=|a_1^{(3)}\wedge a_2^{(1)}|^2=1$. Substituting it into $|W_4|^2+|V_4|^2=1$ in \eqref{holos-thm3-eq5-III} yields $|W_4|^2=0$,
which contradicts $W_4\neq 0$.
Hence we verify $a_{21}^{(1)}=0$.

Now we know $\left\langle a_1^{(4)},~a_2^{(1)}\right\rangle=\left\langle a_2^{(3)},~a_2^{(1)}\right\rangle=0$. The latter gives us $\left\langle a_1^{(3)},~a_1^{(1)}\right\rangle=0$ by $\left\langle W_3,~W_1\right\rangle=0$
in \eqref{holos-thm3-eq5-III}.

Since $V_6=a_{1}^{(3)}\wedge a_2^{(3)}+a_{1}^{(4)}\wedge a_2^{(2)}=\left(a_{21}^{(3)}/a_{11}^{(4)}\cdot a_{1}^{(3)}-a_{2}^{(2)}\right)\wedge a_1^{(4)}=0$,
then there exists $\lambda\in\mathbb{C}$ such that
\begin{equation}\label{III-eq10}
  a_{2}^{(2)}=a_{21}^{(3)}/a_{11}^{(4)}\cdot a_{1}^{(3)}+\lambda a_{1}^{(4)},
\end{equation}
which implies $\left\langle a_1^{(1)},~a_2^{(2)}\right\rangle=0$ by using
$\left\langle a_1^{(1)},~a_1^{(3)}\right\rangle=0=\left\langle a_1^{(1)},~a_1^{(4)}\right\rangle$.
From \eqref{III-eq10}, we have $V_5=a_{1}^{(2)}\wedge a_2^{(3)}+a_{1}^{(3)}\wedge a_2^{(2)}+a_{1}^{(4)}\wedge a_2^{(1)}=\left(a_{21}^{(3)}/a_{11}^{(4)}\cdot a_{1}^{(2)}+\lambda a_{1}^{(3)}-a_{2}^{(1)}\right)\wedge a_1^{(4)}=0$,
then there exists $\mu\in\mathbb{C}$ such that
\begin{equation}\label{III-eq11}
a_{2}^{(1)}=a_{21}^{(3)}/a_{11}^{(4)}\cdot a_{1}^{(2)}+\lambda a_{1}^{(3)}+\mu a_{1}^{(4)}.
\end{equation}
Applying \eqref{III-eq10} and \eqref{III-eq11}, a straightforward calculation shows
\begin{eqnarray*}
V_4&=&a_1^{(1)}\wedge a_2^{(3)}+a_1^{(2)}\wedge a_2^{(2)}+a_1^{(3)}\wedge a_2^{(1)}\\
&=&\left(a_{21}^{(3)}/a_{11}^{(4)}\cdot a_1^{(1)}+\lambda a_1^{(2)}+\mu a_1^{(3)}\right)\wedge a_1^{(4)}.
\end{eqnarray*}
From this, together with $V_2=a_1^{(1)}\wedge a_2^{(1)}$ and $\left\langle a_1^{(4)},~a_1^{(1)}\right\rangle=\left\langle a_1^{(4)},~a_2^{(1)}\right\rangle=0$, we find $\left\langle V_2,~V_4\right\rangle=0$.
Substituting it into
$\left\langle W_2,~W_4\right\rangle+\left\langle V_2,~V_4\right\rangle=0$
in \eqref{holos-thm3-eq5-III} yields
$
\left\langle W_2,~W_4\right\rangle
=0,
$
which leads to $\left\langle a_1^{(2)},~a_1^{(4)}\right\rangle=0$,
i.e.,
$
a_{11}^{(2)}=0.
$

On the other hand, combining $\left\langle V_2,~V_4\right\rangle=0$ and $\left\langle R_2,~R_4\right\rangle=0$ in \eqref{holos-thm3-eq6-III}, we obtain
\begin{equation}\label{III-eq11+1}
\left\langle a_1^{(1)}\wedge a_2^{(1)},~a_1^{(2)}\wedge a_2^{(2)}\right\rangle=0.
\end{equation}
In terms of this, together with $\left\langle a_1^{(1)},~a_2^{(2)}\right\rangle=0$, we get
\begin{equation}\label{III-eq12}
\left\langle a_1^{(1)},~a_1^{(2)}\right\rangle
\left\langle a_2^{(1)},~a_2^{(2)}\right\rangle=0.
\end{equation}
Since $\left\langle W_1,~W_2\right\rangle=0$ in \eqref{holos-thm3-eq5-III}, then
\begin{equation}\label{III-eq13}
\left\langle a_1^{(1)},~a_1^{(2)}\right\rangle
+\left\langle a_2^{(1)},~a_2^{(2)}\right\rangle=0.
\end{equation}
It follows from \eqref{III-eq12} and \eqref{III-eq13} that
\begin{equation}\label{III-eq13+1}
\left\langle a_1^{(1)},~a_1^{(2)}\right\rangle
=\left\langle a_2^{(1)},~a_2^{(2)}\right\rangle=0,
\end{equation}
which implies respectively by \eqref{III-eq10} and \eqref{III-eq11}
\begin{equation}\label{III-eq13+2}
\left\langle a_1^{(1)},~a_2^{(1)}\right\rangle=0
\end{equation}
and
\begin{equation}\label{III-eq14}
\overline{a_{21}^{(3)}/a_{11}^{(4)}}
\left\langle a_2^{(1)},~a_1^{(3)}\right\rangle=0.
\end{equation}
It is worth noting that $\left\langle a_1^{(1)},~a_i^{(\alpha)}\right\rangle=0$ for
$(\alpha,i)\neq(1,1)$.

\textbf{Claim 2:} $a_{21}^{(3)}=0$.

Otherwise if $a_{21}^{(3)}\neq0$, then $\left\langle a_2^{(1)},~a_1^{(3)}\right\rangle=0$
by \eqref{III-eq14}.
Since $S_4=-a_1^{(1)}\wedge a_2^{(1)}\wedge a_1^{(2)},
T_4=-a_1^{(1)}\wedge a_2^{(1)}\wedge a_2^{(2)}$
and $S_6=a_1^{(1)}\wedge a_2^{(1)}\wedge a_1^{(4)}+a_1^{(1)}\wedge a_2^{(3)}\wedge a_1^{(2)}-a_1^{(2)}\wedge a_2^{(1)}\wedge a_1^{(3)},
T_6=a_1^{(3)}\wedge a_2^{(1)}\wedge a_2^{(2)}-a_1^{(1)}\wedge a_2^{(2)}\wedge a_2^{(3)}$,
then a straightforward calculation shows
$$
\left\langle S_4,~S_6\right\rangle
=0=\left\langle T_4,~T_6\right\rangle.
$$
From $\left\langle R_4,~R_6\right\rangle
+\left\langle S_4,~S_6\right\rangle
+\left\langle T_4,~T_6\right\rangle=0$ in \eqref{holos-thm3-eq6-III},
we have $\left\langle R_4,~R_6\right\rangle=0$,
which implies
\\$\left\langle a_1^{(2)}\wedge a_2^{(2)},~a_1^{(3)}\wedge a_2^{(3)}\right\rangle=0$, i.e.,
$$
\left\langle a_1^{(2)},~a_1^{(3)}\right\rangle
a_{21}^{(2)}\overline{a_{21}^{(3)}}=0.
$$
Since $\left\langle R_2,~R_3\right\rangle=0$,
then $\left\langle V_2,~V_3\right\rangle=0$,
which means $\left\langle W_2,~W_3\right\rangle=0$,
i.e.,
$$
\left\langle a_1^{(2)},~a_1^{(3)}\right\rangle
+a_{21}^{(2)}\overline{a_{21}^{(3)}}=0.
$$
So that we obtain
$$
\left\langle a_1^{(2)},~a_1^{(3)}\right\rangle
=0=a_{21}^{(2)}.
$$
Since $\left\langle a_1^{(1)}\wedge a_2^{(2)}+a_1^{(2)}\wedge a_2^{(1)},~a_1^{(2)}\wedge a_2^{(2)}\right\rangle=0$,
then by $\left\langle R_3,~R_4\right\rangle=0$ we have
$\left\langle V_3,~V_4\right\rangle=0$, which implies
$\left\langle W_3,~W_4\right\rangle=0$, i.e., $\left\langle a_1^{(3)},~a_1^{(4)}\right\rangle
=0$. Hence we get $a_{11}^{(3)}=0$.
Using $a_{21}^{(2)}=0=a_{11}^{(3)}$ in \eqref{III-eq10}, we know $\lambda=0$, i.e.,
$$
a_{2}^{(2)}=a_{21}^{(3)}/a_{11}^{(4)}\cdot a_{1}^{(3)}.
$$
Then by $a_{21}^{(1)}=0=a_{11}^{(2)}$ in \eqref{III-eq11}, we have $\mu=0$, i.e.,
$$
a_{2}^{(1)}=a_{21}^{(3)}/a_{11}^{(4)}\cdot a_{1}^{(2)}.
$$
Since $V_4=a_1^{(1)}\wedge a_2^{(3)}$ and
$$
R_5=-2a_1^{(4)}\wedge a_2^{(1)},
S_5=0,
T_5=-2a_1^{(1)}\wedge a_2^{(1)}\wedge a_2^{(3)},
$$
then using $|W_4|^2+|V_4|^2=1$ in \eqref{holos-thm3-eq5-III} and
$|R_5|^2+|S_5|^2+|T_5|^2=4c$ in \eqref{holos-thm3-eq6-III}, we obtain
\begin{equation*}
  |a_{11}^{(4)}|^2+|a_{1}^{(1)}|^2|a_{21}^{(3)}|^2=1,
  ~|a_{11}^{(4)}|^2|a_{2}^{(1)}|^2
+|a_{1}^{(1)}|^2|a_{2}^{(1)}|^2|a_{21}^{(3)}|^2=c,
\end{equation*}
which implies $|a_{2}^{(1)}|^2=c$.
Substituting this into $|R_2|^2=|a_{1}^{(1)}|^2|a_{2}^{(1)}|^2=c$ in \eqref{holos-thm3-eq6-III} yields
$|a_{1}^{(1)}|^2=1$.
From $|W_1|^2=|a_1^{(1)}|^2+|a_{2}^{(1)}|^2=4$ in \eqref{holos-thm3-eq5-III},
we have $|a_{2}^{(1)}|^2=3$, which implies $c=3$.
It follows that $|V_2|^2=|R_2|^2=3$, which leads to
$|W_2|^2=|a_1^{(2)}|^2+|a_{2}^{(2)}|^2=3$ by $|W_2|^2+|V_2|^2=6$ in \eqref{holos-thm3-eq5-III}.
On the other hand, that $|R_3|^2=4|a_{1}^{(1)}|^2|a_{2}^{(2)}|^2=12$ in \eqref{holos-thm3-eq6-III}
and $|a_{1}^{(1)}|^2=1$ tells us $|a_{2}^{(2)}|^2=3$, which shows $|a_1^{(2)}|^2=0$, i.e., $a_2^{(1)}=0$.
It contradicts $|a_{2}^{(1)}|^2=3$. Thus we verify $a_{21}^{(3)}=0$.

At this time, equations \eqref{III-eq10} and \eqref{III-eq11} reduce to
$$
a_{2}^{(2)}=a_{21}^{(2)}/a_{11}^{(4)}\cdot a_{1}^{(4)},
~a_{2}^{(1)}=a_{21}^{(2)}/a_{11}^{(4)}\cdot a_{1}^{(3)}+\mu a_{1}^{(4)}.
$$
Since $V_3=a_1^{(1)}\wedge a_2^{(2)}+a_1^{(2)}\wedge a_2^{(1)},~V_4=a_1^{(2)}\wedge a_2^{(2)}+a_1^{(3)}\wedge a_2^{(1)},$ then we observe that $\left\langle V_3,~V_4\right\rangle=0$, which implies
$\left\langle W_3,~W_4\right\rangle=0$, i.e., $a_{11}^{(3)}=0$. So that $a_{2}^{(1)}=a_{21}^{(2)}/a_{11}^{(4)}\cdot a_{1}^{(3)}$.

Since $V_2=R_2=a_1^{(1)}\wedge a_2^{(1)}$, then using $|V_2|^2=|R_2|^2=c$ in
\eqref{holos-thm3-eq6-III} yields
\begin{equation}\label{III-eq15}
  |a_{1}^{(1)}|^2|a_{2}^{(1)}|^2=c.
\end{equation}
From $|W_2|^2+|V_2|^2=6$ in \eqref{holos-thm3-eq5-III}, we know that $|W_2|^2=6-c$, i.e.,
\begin{equation}\label{III-eq16}
  |a_{1}^{(2)}|^2+|a_{21}^{(2)}|^2=6-c.
\end{equation}
Because $R_4=4a_1^{(2)}\wedge a_2^{(2)},S_4=-a_1^{(1)}\wedge a_2^{(1)}\wedge a_1^{(2)},
T_4=-a_1^{(1)}\wedge a_2^{(1)}\wedge a_2^{(2)},$
then applying $|R_4|^2+|S_4|^2+|T_4|^2=6c$ in \eqref{holos-thm3-eq6-III} gives us
\begin{equation}\label{III-eq17}
  16|a_{1}^{(2)}|^2|a_{21}^{(2)}|^2+|a_{1}^{(1)}|^2|a_{2}^{(1)}|^2
|a_{1}^{(2)}|^2+|a_{1}^{(1)}|^2|a_{2}^{(1)}|^2|a_{21}^{(2)}|^2=6c.
\end{equation}
Substituting \eqref{III-eq15} and \eqref{III-eq16} into \eqref{III-eq17} leads to
\begin{equation}\label{III-eq18}
  |a_{1}^{(2)}|^2|a_{21}^{(2)}|^2=\frac{c^2}{16},
\end{equation}
which implies $|V_4|^2=\frac{c^2}{16}$. From this, together with $|W_4|^2+|V_4|^2=1$ in \eqref{holos-thm3-eq5-III}, we get
\begin{equation}\label{III-eq19}
  |a_{11}^{(4)}|^2=1-\frac{c^2}{16},
\end{equation}
which means that $0<c<4$ by $a_{11}^{(4)}>0$.
On the other hand, since $R_5=-2a_1^{(4)}\wedge a_2^{(1)},S_5=2a_1^{(1)}\wedge a_1^{(2)}\wedge a_2^{(2)},
T_5=2a_2^{(1)}\wedge a_1^{(2)}\wedge a_2^{(2)}$, then it follows from $|R_5|^2+|S_5|^2+|T_5|^2=4c$ in
\eqref{holos-thm3-eq6-III} that
\begin{equation}\label{III-eq20}
  |a_{11}^{(4)}|^2|a_{2}^{(1)}|^2
+|a_{1}^{(1)}|^2|a_{1}^{(2)}|^2|a_{21}^{(2)}|^2
+|a_{2}^{(1)}|^2|a_{1}^{(2)}|^2|a_{21}^{(2)}|^2=c.
\end{equation}
Now substituting \eqref{III-eq18}, \eqref{III-eq19} and $|W_1|^2=|a_{1}^{(1)}|^2+|a_{2}^{(1)}|^2=4$ into
\eqref{III-eq20} yields
\begin{equation}\label{III-eq21}
  |a_{2}^{(1)}|^2=\frac{4c}{c+4},
\end{equation}
which implies $|a_{1}^{(1)}|^2=\frac{16}{c+4}$.
Then it follows from \eqref{III-eq15} that $\frac{64c}{(c+4)^2}=c$, which contradicts that $0<c<4$.
Thus there does not exist the case of $W_4\neq 0$.
\end{proof}

\begin{proof}{\emph{Proof of Theorem \ref{main-theorem2}}}
It is enough to consider the irreducible case.
Assume that $\varphi$ generates the following harmonic sequence (cf. \cite{Burstall-Wood 1986},\cite{Chern-Wolfson 1987})
\begin{equation}
0\stackrel{A_{\varphi_0}''}{\longleftarrow} \underline{\varphi}_0 = \underline{\varphi}
\stackrel{A_{\varphi_0}'}{\longrightarrow} \underline{\varphi}_1
\stackrel{A_{\varphi_1}'}{\longrightarrow} \underline{\varphi}_2
\stackrel{A_{\varphi_{2}}'}{\longrightarrow} \underline{\varphi}_{3}
\stackrel{A_{\varphi_3}'}{\longrightarrow} 0,
\label{t2-eq1}
\end{equation}
where $\text{rank}(\underline{\varphi}_0)=\text{rank}(\underline{\varphi}_1)=2$ and $\left(\text{rank}(\underline{\varphi}_2),\text{rank}(\underline{\varphi}_3)\right)=(1,1)
~\text{or}~(2,0)$.
We will prove $\mathrm{d}=4$ in the above two cases.

If $\left(\text{rank}(\underline{\varphi}_2),\text{rank}(\underline{\varphi}_3)\right)=(1,1)$,
then $\underline{\varphi}_2=\underline{f}_{m-1}^{(m)},~\underline{\varphi}_3=\underline{f}_{m}^{(m)}$
belong the following harmonic sequence in $\mathbb{C}P^m$
\begin{equation}
0\stackrel{A_{f_0^{(m)}}''}{\longleftarrow} \underline{f}_0^{(m)}
\stackrel{A_{f_0^{(m)}}'}{\longrightarrow} \underline{f}_1^{(m)}
\stackrel{A_{f_1^{(m)}}'}{\longrightarrow} \cdots
\stackrel{A_{f_{m-2}^{(m)}}'}{\longrightarrow} \underline{f}_{m-1}^{(m)}
\stackrel{A_{f_{m-1}^{(m)}}'}{\longrightarrow} \underline{f}_{m}^{(m)}
\stackrel{A_{f_m^{(m)}}'}{\longrightarrow} 0
\label{t2-eq2}
\end{equation}
for $m=5,4,3$.
Notice that Fei-He \cite{Fei-He 2019} discuss this case under the assumption that $\varphi$ is unramified.
Since the condition of constant square norm of the second fundamental form implies the unramified one, then we give an explicit proof.

When $m=5$, we have a diagram below about the harmonic sequence in $\mathbb{C}P^{5}$ (cf. \cite{Burstall-Wood 1986}),
\begin{eqnarray*}
\xymatrix{
  \underline{e}_2 \ar[dr] \ar[r]  & \underline{e}_4=\underline{f}_3^{(5)} \ar[r] & \underline{\varphi}_2=\underline{f}_4^{(5)} \ar[r]& \underline{\varphi}_3=\underline{f}_5^{(5)} \ar[r] &0.\\
  \underline{e}_1 \ar[u] \ar[r]_{} & \underline{e}_3  \ar[u]   &    &    &   &   }
\end{eqnarray*}
It follows that
$$
v_1\wedge v_2\wedge \partial_z v_1\wedge\partial_z v_2
=h(z)f_0^{(5)}\wedge\partial_zf_0^{(5)}\wedge\partial^2_zf_0^{(5)}\wedge \partial^3_zf_0^{(5)},
$$
where $h(z)$ is a polynomial. Notice that $f_0^{(5)}\wedge\partial_zf_0^{(5)}\wedge\partial^2_zf_0^{(5)}\wedge \partial^3_zf_0^{(5)}$ is the $3$-osculating curve $\sigma_3$ of $f_0^{(5)}$. It follows from \eqref{eq2+} that $h$ is a non-zero constant and $\sigma_3$ has constant curvature. Using Y.B. Shen's result
(\cite{Shen 1996}, Theorem 6.2), we know the harmonic sequence \eqref{t2-eq2} in $\mathbb{C}P^{5}$ is the Veronese sequence, up to $U(6)$.
Then we obtain by \eqref{eq2+} and (3.24) in \cite{Bolton-Jensen-Rigoli-Woodward 1988}
$$
2\mathrm{d}-4=\delta_3^{(5)}=8,
$$
where $\delta_3^{(5)}$ is the degree of $\sigma_3$.
It follows that $\mathrm{d}=6$.
In this case, there exist polynomials $P_1,P_2,P_3$ such that
$$
v_1\wedge v_2=P_1f_0^{(5)}\wedge\partial_zf_0^{(5)}+P_2f_0^{(5)}\wedge\partial^2_zf_0^{(5)}+P_3\partial_zf_0^{(5)}\wedge\partial^2_zf_0^{(5)}.
$$
We choose a rectangular coordinate system in $\mathbb{C}^6$ such that
$$
f_0^{(5)}=\epsilon_1+\sqrt{5}z\epsilon_2+\sqrt{10}z^2\epsilon_3+\sqrt{10}z^3\epsilon_4+\sqrt{5}z^4\epsilon_5+z^5\epsilon_6.
$$
A straightforward calculation shows that the degree of $f_0^{(5)}\wedge\partial_zf_0^{(5)}$,
$f_0^{(5)}\wedge\partial^2_zf_0^{(5)}$ and
$\partial_zf_0^{(5)}\wedge\partial^2_zf_0^{(5)}$
is $8$, $7$ and $6$ respectively, and
\begin{eqnarray*}
f_0^{(5)}\wedge\partial_zf_0^{(5)}
&=&\sqrt{5}\epsilon_1\wedge\epsilon_2+2\sqrt{10}z\epsilon_1\wedge\epsilon_3+3\sqrt{10}z^2\epsilon_1\wedge\epsilon_4
+4\sqrt{5}z^3\epsilon_1\wedge\epsilon_5+5z^4\epsilon_1\wedge\epsilon_6\\
&&+\sqrt{50}z^2\epsilon_2\wedge\epsilon_3+2\sqrt{50}z^3\epsilon_2\wedge\epsilon_4
+15z^4\epsilon_2\wedge\epsilon_5+4\sqrt{5}z^5\epsilon_2\wedge\epsilon_6\\
&&+10z^4\epsilon_3\wedge\epsilon_4
+2\sqrt{50}z^5\epsilon_3\wedge\epsilon_5+3\sqrt{10}z^6\epsilon_3\wedge\epsilon_6\\
&&+\sqrt{50}z^6\epsilon_4\wedge\epsilon_5+2\sqrt{10}z^7\epsilon_4\wedge\epsilon_6+\sqrt{5}z^8\epsilon_5\wedge\epsilon_6,
\end{eqnarray*}
\begin{eqnarray*}
f_0^{(5)}\wedge\partial^2_zf_0^{(5)}
&=&2\sqrt{10}\epsilon_1\wedge\epsilon_3+6\sqrt{10}z\epsilon_1\wedge\epsilon_4
+12\sqrt{5}z^2\epsilon_1\wedge\epsilon_5+20z^3\epsilon_1\wedge\epsilon_6\\
&&+2\sqrt{50}z\epsilon_2\wedge\epsilon_3+6\sqrt{50}z^2\epsilon_2\wedge\epsilon_4
+60z^3\epsilon_2\wedge\epsilon_5+20\sqrt{5}z^4\epsilon_2\wedge\epsilon_6\\
&&+40z^3\epsilon_3\wedge\epsilon_4
+10\sqrt{50}z^4\epsilon_3\wedge\epsilon_5+18\sqrt{10}z^5\epsilon_3\wedge\epsilon_6\\
&&+6\sqrt{50}z^5\epsilon_4\wedge\epsilon_5+14\sqrt{10}z^6\epsilon_4\wedge\epsilon_6+8\sqrt{5}z^7\epsilon_5\wedge\epsilon_6,
\end{eqnarray*}
\begin{eqnarray*}
\partial_zf_0^{(5)}\wedge\partial^2_zf_0^{(5)}
&=&2\sqrt{50}\epsilon_2\wedge\epsilon_3+6\sqrt{50}z\epsilon_2\wedge\epsilon_4
+60z^2\epsilon_2\wedge\epsilon_5+20\sqrt{5}z^3\epsilon_2\wedge\epsilon_6\\
&&+60z^2\epsilon_3\wedge\epsilon_4
+16\sqrt{50}z^3\epsilon_3\wedge\epsilon_5+30\sqrt{10}z^4\epsilon_3\wedge\epsilon_6\\
&&+12\sqrt{50}z^4\epsilon_4\wedge\epsilon_5+30\sqrt{10}z^5\epsilon_4\wedge\epsilon_6+20\sqrt{5}z^6\epsilon_5\wedge\epsilon_6.
\end{eqnarray*}
Since $[v_1\wedge v_2]$ is a constantly curved holomorphic curve of degree $6$, then under the standard 
basis $\epsilon_i\wedge\epsilon_j,~1\leq i<j\leq 6$ in the
lexicographic order, we can write $v_1\wedge v_2=PHV_0^{(6)}$,
where $P$ is a polynomial, $V_0^{(6)}=(1,\sqrt{6}z,\sqrt{15}z^2,\sqrt{20}z^3,\sqrt{15}z^4,\sqrt{6}z^5,z^6)^{T}$
and $H=(H_{\alpha j})_{\alpha=0,\cdots,14,j=0,\cdots,6}$ is a $15\times 7$ matrix whose column vectors are mutually orthonormal.
A straightforward calculation shows
$$
\sqrt{5}P_1=P(H_{00}+\cdots+H_{06}z^6),
~2\sqrt{10}zP_1+2\sqrt{10}P_2=P(H_{10}+\cdots+H_{16}z^6),
$$ 
$$
\sqrt{50}z^2P_1+2\sqrt{50}zP_2+2\sqrt{50}P_3=P(H_{50}+\cdots+H_{56}z^6),
$$ 
which implies $P| P_j(j=1,2,3)$. Set $P_j=P\tilde{P}_j$, then we have
$$
\tilde{P}_1f_0^{(5)}\wedge\partial_zf_0^{(5)}
+\tilde{P}_2f_0^{(5)}\wedge\partial^2_zf_0^{(5)}
+\tilde{P}_3\partial_zf_0^{(5)}\wedge\partial^2_zf_0^{(5)}
=HV_0^{(6)}.
$$
It follows from the above formulas that
the degree of $\tilde{P}_1,\tilde{P}_2,\tilde{P}_3$ is $0,1,2$ respectively.
Set $\tilde{P}_1=\tilde{p}_{10}\neq 0$, $\tilde{P}_{2}=\tilde{p}_{20}+\tilde{p}_{21}z$ 
and $\tilde{P}_3=\tilde{p}_{30}+\tilde{p}_{31}z+\tilde{p}_{32}z^2$,
then that the terms of $z^7,z^8$ in $v_1\wedge v_2$ are vanishing tells us
\begin{equation}\label{eq+1}
\tilde{p}_{10}+7\tilde{p}_{21}+15\tilde{p}_{32}=0,
\end{equation}
\begin{equation}\label{eq+2}
\tilde{p}_{10}+8\tilde{p}_{21}+20\tilde{p}_{32}=0,
\end{equation}
and
\begin{equation}\label{eq+3}
2\tilde{p}_{20}+5\tilde{p}_{31}=0.
\end{equation}
From \eqref{eq+1} and \eqref{eq+2}, we get
\begin{equation}\label{eq+4}
\tilde{p}_{21}=-\frac{1}{4}\tilde{p}_{10},~\tilde{p}_{32}=\frac{1}{20}\tilde{p}_{10}.
\end{equation}
At this time $H$ is given by
$$
H=\begin{pmatrix}
H_{00} & 0 & 0 & 0 & 0 & 0 & 0\\
H_{10} & H_{11} & 0 & 0 & 0 & 0 & 0\\
0 & H_{21} &  H_{22} & 0 & 0 & 0 & 0\\
0 & 0 &  H_{32} & H_{33} & 0 & 0 & 0\\
0 & 0 &  0 &  H_{43} & H_{44} & 0 & 0\\
H_{50} & H_{51} & H_{52} & 0 & 0 & 0 & 0\\
0 & H_{61} & H_{62} & H_{63} & 0 & 0 & 0\\
0 & 0 & H_{72} & H_{73} & H_{74} & 0 & 0\\
0 & 0 & 0 & H_{83} & H_{84} & H_{85} & 0\\
0 & 0 & H_{92} & H_{93} & H_{94} & 0 & 0\\
0 & 0 & 0 & H_{10,3} & H_{10,4} & H_{10,5} & 0\\
0 & 0 & 0 & 0 & H_{11,4} & H_{11,5} & H_{11,6}\\
0 & 0 & 0 & 0 & H_{12,4} & H_{12,5} & H_{12,6}\\
0 & 0 & 0 & 0 & 0 & H_{13,5} & H_{13,6}\\
0 & 0 & 0 & 0 & 0 & 0 & H_{14,6}\\
\end{pmatrix}.
$$
Using $\eqref{eq+4}$, we know that $H_{11}H_{52}\neq 0$.
Since the column vectors of $H$
are mutually orthonormal, then $\overline{H_{50}}H_{52}=0$, which implies
$H_{50}=0$, i.e., $\tilde{p}_{30}=0$.
From $\overline{H_{10}}H_{11}=0$, we have $H_{10}=0$, i.e., $\tilde{p}_{20}=0$.
Substituting this into \eqref{eq+3} yields $\tilde{p}_{31}=0$.
Now $H$ becomes
$$
H=\begin{pmatrix}
H_{00} & 0 & 0 & 0 & 0 & 0 & 0\\
0 & H_{11} & 0 & 0 & 0 & 0 & 0\\
0 & 0 &  H_{22} & 0 & 0 & 0 & 0\\
0 & 0 &  0 & H_{33} & 0 & 0 & 0\\
0 & 0 &  0 &  0 & H_{44} & 0 & 0\\
0 & 0 & H_{52} & 0 & 0 & 0 & 0\\
0 & 0 & 0 & H_{63} & 0 & 0 & 0\\
0 & 0 & 0 & 0 & H_{74} & 0 & 0\\
0 & 0 & 0 & 0 & 0 & H_{85} & 0\\
0 & 0 & 0 & 0 & H_{94} & 0 & 0\\
0 & 0 & 0 & 0 & 0 & H_{10,5} & 0\\
0 & 0 & 0 & 0 & 0 & 0 & H_{11,6}\\
0 & 0 & 0 & 0 & 0 & 0 & H_{12,6}\\
0 & 0 & 0 & 0 & 0 &0 & 0\\
0 & 0 & 0 & 0 & 0 & 0 & 0\\
\end{pmatrix}.
$$
Since $H_{00}=\sqrt{5}\tilde{p}_{10}$ and $H_{11}=\frac{3\sqrt{10}}{2}\tilde{p}_{10}$,
then it follows from $|H_{00}|^2=1$ that $|\tilde{p}_{10}|^2=\frac{1}{5}$, which implies that
$|H_{11}|^2=\frac{3}{2}$. It contradicts that $|H_{11}|^2=1$.
Thus there does not exist the case of $m=5$.

When $m=4$, the trivial bundle $S^2\times \mathbb{C}^6$ over $S^2$ has a decomposition $S^2\times \mathbb{C}^6=
S^2\times \mathbb{C}^{5}\oplus S^2\times \mathbb{C}$. In this case, the diagram reduces to
\begin{eqnarray*}
\xymatrix{
  \underline{e}_2 \ar[dr] \ar[r]  & \underline{e}_4=\underline{f}_2^{(4)} \ar[r] & \underline{\varphi}_2=\underline{f}_3^{(4)} \ar[r]&\underline{\varphi}_3=\underline{f}_4^{(4)} \ar[r] &0,\\
  \underline{e}_1 \ar[u] \ar[r]_{} & \underline{e}_3  \ar[u]   &    &    &   &   }
\end{eqnarray*}
where $\mathrm{span}_\mathbb{C}\left\{e_1,e_2,e_3\right\}
=\mathrm{span}_\mathbb{C}\left\{v_0,f_0^{(4)},f_1^{(4)}\right\}$
with $v_0=(0,0,0,0,0,1)$.
Similarly we get
$$
2\mathrm{d}-4=\delta_2^{(4)}=6,
$$
which implies $\mathrm{d}=5$.
In this case, there exist polynomials $P_1,P_2,P_3$ such that
$$
v_1\wedge v_2=P_1f_0^{(4)}\wedge v_0+P_2\partial_zf_0^{(4)}\wedge v_0+P_3f_0^{(4)}\wedge\partial_zf_0^{(4)}.
$$
We choose a rectangular coordinate system in $\mathbb{C}^5$ such that
$$
f_0^{(4)}=\epsilon_1+2z\epsilon_2+\sqrt{6}z^2\epsilon_3+2z^3\epsilon_4+z^4\epsilon_5.
$$
A straightforward calculation shows that the degree of $f_0^{(4)}\wedge v_0$,
$\partial_zf_0^{(4)}\wedge v_0$ and
$f_0^{(4)}\wedge\partial_zf_0^{(4)}$
is $4$, $3$ and $6$ respectively, and
\begin{eqnarray*}
f_0^{(4)}\wedge v_0
=\epsilon_1\wedge\epsilon_6+2z\epsilon_2\wedge\epsilon_6+\sqrt{6}z^2\epsilon_3\wedge\epsilon_6
+2z^3\epsilon_4\wedge\epsilon_6+z^4\epsilon_5\wedge\epsilon_6,
\end{eqnarray*}
\begin{eqnarray*}
\partial_zf_0^{(4)}\wedge v_0
=2\epsilon_2\wedge\epsilon_6+2\sqrt{6}z\epsilon_3\wedge\epsilon_6
+6z^2\epsilon_4\wedge\epsilon_6+4z^3\epsilon_5\wedge\epsilon_6,
\end{eqnarray*}
\begin{eqnarray*}
f_0^{(4)}\wedge\partial_zf_0^{(4)}
&=&2\epsilon_1\wedge\epsilon_2+2\sqrt{6}z\epsilon_1\wedge\epsilon_3
+6z^2\epsilon_1\wedge\epsilon_4+4z^3\epsilon_1\wedge\epsilon_5\\
&&+2\sqrt{6}z^2\epsilon_2\wedge\epsilon_3
+8z^3\epsilon_2\wedge\epsilon_4+6z^4\epsilon_2\wedge\epsilon_5\\
&&+2\sqrt{6}z^4\epsilon_3\wedge\epsilon_4+2\sqrt{6}z^5\epsilon_3\wedge\epsilon_5+2z^6\epsilon_4\wedge\epsilon_5.
\end{eqnarray*}
Since $[v_1\wedge v_2]$ is a constantly curved holomorphic curve of degree $5$, then under the standard
basis $\epsilon_i\wedge\epsilon_j,~1\leq i<j\leq 6$ in the
lexicographic order, we can write $v_1\wedge v_2=PHV_0^{(5)}$,
where $P$ is a polynomial, $V_0^{(5)}=(1,\sqrt{5}z,\sqrt{10}z^2,\sqrt{10}z^3,\sqrt{5}z^4,z^5)^{T}$
and $H=(H_{\alpha j})_{\alpha=0,\cdots,14,j=0,\cdots,5}$ is a $15\times 6$ matrix whose column vectors are mutually orthonormal.
A similar discussion shows that $P| P_j(j=1,2,3)$. Set $P_j=P\tilde{P}_j$, then we have
$$
\tilde{P}_1f_0^{(4)}\wedge v_0
+\tilde{P}_2\partial_zf_0^{(4)}\wedge v_0
+\tilde{P}_3f_0^{(4)}\wedge\partial_zf_0^{(4)}
=HV_0^{(5)}.
$$
It follows from the above formulas that
$\tilde{P}_3=0$.
So $\underline{\varphi}=\text{span}\left\{v_0,\tilde{P}_1f_0^{(4)}+\tilde{P}_2\partial_zf_0^{(4)}\right\}$, 
which implies that $\varphi$ is reducible.
Hence there does not exist the case of $m=4$.

When $m=3$, then the trivial bundle $S^2\times \mathbb{C}^6$ over $S^2$ has a decomposition $S^2\times \mathbb{C}^6=
S^2\times \mathbb{C}^{4}\oplus S^2\times \mathbb{C}^2$. In this case,
we have
$$
2\mathrm{d}-4=\delta_1^{(3)}=4,
$$
which implies $\mathrm{d}=4$.

If $\left(\text{rank}(\underline{\varphi}_2),\text{rank}(\underline{\varphi}_3)\right)=(2,0)$,
then we choose a local unitary frame $e=\left\{e_1,\cdots,e_6\right\}$
along $\varphi$ so that
$$
\underline{\varphi}_0=\text{span}\{e_{1},e_{2}\},
~\underline{\varphi}_1=\text{span}\{e_{3},e_{4}\},
~\underline{\varphi}_2=\text{span}\{e_{5},e_{6}\}.
$$
Under such frame, the pull back of (right invariant) Maurer-Cartan forms which are denoted by $\omega=(\omega_{AB})$ are
\begin{equation*}
\left(
  \begin{array}{ccccccccc}
    \Omega_1 & A_1\phi & {} \\
    -A_1^*\bar{\phi} & \Omega_2 & A_2\phi \\
    {} & -A_2^*\bar{\phi} & \Omega_3 \\
  \end{array}
\right).
\end{equation*}
Notice that the unitary frame we choose is determined up to a transformation of
the group $U(2)\times U(2)\times
U(2)$, so $|\det A_i|$ $(i=1,2)$ are global invariants of analytic type on $S^2$ vanishing only at isolated points, and away from their zeros, they satisfy (cf. \cite{Chi-Zheng 1989}, \cite{Fei-Jiao-Xu 2011})
\begin{equation}\label{t2-eq3}
\Delta\log|\det A_1|=2K-4L_1+2L_{2},
\end{equation}
\begin{equation}\label{t2-eq4}
\Delta\log|\det A_2|=2K+2L_1-4L_{2},
\end{equation}
where $K=4/\mathrm{d}$, $L_i=trA_iA_i^*~(i=1,2)$ are also globally defined invariants on $S^2$ with $L_1=1$, and $\Delta$ is Laplace-Beltrami operator with respect to $\mathrm{d}\cdot ds_{S^2}^2$.
It follows from \eqref{t2-eq3} and \eqref{t2-eq4} that
\begin{equation}\label{t2-eq5}
\Delta\log|\det A_1|^2|\det A_2|=\frac{6(4-\mathrm{d})}{\mathrm{d}}.
\end{equation}
Since $\varphi$ is totally unramified, then $|\det A_1|^2|\det A_2|$ has no zeros on $S^2$.
Applying the maximum principle to \eqref{t2-eq5},
we get $\mathrm{d}=4$.
The conclusion follows from Theorem \ref{main-theorem1}.
\end{proof}

\begin{remark}
In the irreducible case of Theorem \ref{main-theorem2}, the condition of constant square norm of the second fundamental form is necessary. Jiao \cite{Jiao 2008} gave an example of irreducible constantly curved holomorphic two-sphere of $\mathrm{d}=4$ in $G(2,6)$, as follows,
\begin{equation}
\varphi=\begin{bmatrix}
1&0&\frac{1}{\sqrt{2}}z&\frac{\sqrt{31}}{2\sqrt{7}}z^2&\frac{9}{2\sqrt{7}}z^2
&0\\
0&1&0&0&\frac{\sqrt{7}}{\sqrt{2}}z&\frac{1}{2}z^2
\end{bmatrix}
\end{equation}
with
$$
|\det A_1|^2=\frac{112+1024z\bar{z}+1176z^2\bar{z}^2
+376z^3\bar{z}^3+31z^4\bar{z}^4}{1024(1+z\bar{z})^4},
$$
which is not a constant.
\end{remark}

\end{document}